\documentclass{amsproc}
\usepackage{amsmath}
\usepackage{amscd}
\usepackage{amssymb}

\usepackage{amsthm}
    
\usepackage[usenames]{color}    
\newtheorem{Def}{Definition}[section]     
\newtheorem{Prop}[Def]{Proposition}
\newtheorem{Lemma}[Def]{Lemma}     
\newtheorem{Thm}[Def]{Theorem} 
 
\newtheorem{Rem}[Def]{Remark}

\newcommand{\Q}{\mathbb{Q}}
\newcommand{\C}{\mathbb{C}}
\newcommand{\R}{\mathbb{R}}
\newcommand{\Z}{\mathbb{Z}}

\newcommand{\Ee}{\tilde{\mathbb E}}
\newcommand{\Hh}{\mathbb{H}}

\newcommand{\EW}{\Ee//\tilde{W}}

\newcommand{\Der}{\mathrm{Der}}
\newcommand{\Specan}{\mathrm{Specan}}

\newcommand{\p}{\partial}

\newcommand{\e}{\epsilon}

\newcommand{\s}{\sigma}

\newcommand{\n}{\mathfrak {n}}

\newcommand{\de}{{\delta}}

\newcommand{\Oo}{\mathcal {O}}

\newcommand{\Ff}{\mathcal {F}}
\newcommand{\Ss}{\hbox{\boldmath$\mathcal {S}$}}

\newcommand{\bp}{\begin{pmatrix}}
\newcommand{\ep}{\end{pmatrix}}

\numberwithin{equation}{section}

\newcounter{CounterEQUlabel}
\newcommand{\EQUlabel}[1]{\label{#1}
	\ifcase \theCounterEQUlabel
		\relax
	\or
		\hspace{1em}\mbox{\tiny$\langle$\rmfamily#1$\rangle$}
		\index{zzz#1@#1}
	\fi }	
	\newcounter{CounterEQUref}
	\newcounter{CounterEQUpageref}
	\newcommand{\EQUref}[1]{
		\ifcase \theCounterEQUref     \relax   \or {\tiny[#1]}\,\fi
		\ifcase \theCounterEQUpageref (\ref{#1})\or (\ref{#1}\,(p.\pageref{#1})) \fi}

\title[Frobenius manifolds for elliptic root systems]
{Frobenius manifolds for elliptic root systems}
\author{Ikuo Satake}
\address{Department of Mathematics, Osaka University, 
Toyonaka, Osaka, 560-0043, JAPAN}
\date{}
\begin{document}
\maketitle
%
\begin{abstract}
In this paper, we show that 
the quotient space of the domain 
by the reflection group 
for an elliptic root system 
has a structure of Frobenius manifold 
for the case of codimension 1. 
We also give a characterization of this Frobenius manifold 
structure under some suitable condition. 
\end{abstract}
%
%

\section{Introduction}
A Frobenius manifold is a complex manifold 
whose holomorphic tangent bundle has the following 
structures: a flat holomorphic metric 
and a product satisfying some integrable condition. 
It is introduced by Dubrovin in order to represent 
the integrable structures of the topological field theory. 
Its construction is important from the viewpoint of mirror symmetry 
(cf. \cite{Givental}). 

For a Frobenius manifold, the notion of 
``the intersection form" is defined. 
It is a holomorphic symmetric tensor 
on the cotangent bundle of the Frobenius manifold. 
This is introduced in \cite{Saito} in the context of singularity 
theory and defined in \cite{Dubrovin} for Frobenius manifolds. 
Then we set the following problem:

Problem: {\it Let $(M,I^*)$ be a suitable pair where $M$ is a complex manifold
and $I^*$ is a holomorphic symmetric tensor 
on the cotangent bundle of $M$. 
Construct the Frobenius manifold structure on $M$ such that 
its intersection form coincides with the tensor $I^*$.}

If $M$ is the complex orbit space of a finite irreducible Coxeter group 
and $I^*$ is the tensor descended from the standard holomorphic metric, 
this problem is solved by Saito \cite{finite}, Dubrovin \cite{Dubrovin}.

In this paper, we solve the problem for 
the complex orbit space of the reflection group for 
an elliptic root system of codimension 1 with the tensor $I^*$ 
descended from the standard holomorphic metric. 
It is a natural generalization of the complex 
orbit space of a finite Coxeter group. 
We also show a strong uniqueness theorem 
which asserts that the structure of Frobenius manifold 
whose intersection form coincides with $I^*$ is unique 
under the condition that only the Euler vector field is fixed.

We remark that the complex orbit space of the reflection group for an 
elliptic root system appears in various contexts and 
is studied from various viewpoints, for example, 
invariant theory 
(I. N. Bern\v ste\u\i n and O. V. \v Svarcman \cite{Chevalley3}, 
\cite{Chevalley4}, 
E. Looijenga \cite{Chevalley1}, 
K. Saito \cite{extendedII}), 
characters of an affine Lie algebra 
(V. G. Kac and D. H. Peterson \cite{Chevalley2}, 
P. Slodowy \cite{Slodowy}), 
the adjoint quotient of an elliptic Lie group
(S. Helmke and P. Slodowy \cite{Helmke}), 
Jacobi forms (K. Wirthm\"uller \cite{Chevalley5}) and 
moduli spaces of the $G$-principal bundles over  
elliptic curves 
(R. Friedman and J.W. Morgan \cite{Chevalley6}). 

We shall explain the outline of our method of 
the construction of the Frobenius manifold, 
which is almost parallel to the finite Coxeter group case \cite{Dubrovin}. 
Among the structures of the Frobenius manifold, 
the flat metric $J$ is already constructed by Saito \cite{extendedII}. 
A multiplication on the tangent bundle is constructed as follows. 

If we have a Frobenius manifold in general, 
we have the following diagram:
$$
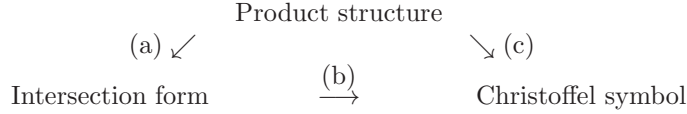
\begin{array}{ccccc}
& \mbox{Product structure} & \\
\qquad\qquad\hbox{(a)}\swarrow && 
\searrow \hbox{(c)}\qquad\qquad\qquad \\
\mbox{Intersection form} & \stackrel{\hbox{(b)}}{\longrightarrow}
& \mbox{Christoffel symbol}
\end{array}.
$$
This diagram means that the product structure gives the intersection 
form (step (a)), the intersection form defines its Levi-Civita connection 
and gives the Christoffel symbols (step (b)), 
therefore the product structure gives the Christoffel symbols (step (c)). 

For step (c), the structure coefficients $C^{\alpha \beta}_{\gamma}$ 
of the product and the Christoffel symbols 
$\Gamma^{\alpha \beta}_{\gamma}$ have the simple relation 
(\cite[p.194, Lemma 3.4]{Dubrovin}):
\begin{equation}\EQUlabel{1.100.100}
\Gamma^{\alpha \beta}_{\gamma}
=
(d^{\beta}+\frac{1-D}{2})C^{\alpha \beta}_{\gamma}
\end{equation}
with respect to the flat coordinates of the Frobenius manifold, 
where $D$ is a degree of the flat metric (see Definition \ref{3.333})
of the Frobenius manifold 
and $d^{\beta}$ is a degree of the homogeneous 
flat coordinate $t^{\beta}$. 

Thanks to the equation (\ref{1.100.100}), the converse construction 
of the step (c) is possible for some cases, which gives a clue to solve 
our problem. 

For the case of the complex orbit space of a finite irreducible 
Coxeter group, the coefficients $d^{\beta}+\frac{1-D}{2}$ 
in (\ref{1.100.100})
are all non-zero. Thus we could define $C^{\alpha \beta}_{\gamma}$ 
by $\Gamma^{\alpha \beta}_{\gamma}$. 
Furthermore we see that $C^{\alpha \beta}_{\gamma}$ satisfies 
the conditions of Frobenius manifold 
by the method of flat pencils.

For the case of the complex orbit space of the reflection 
group for an elliptic root system of codimension 1, 
some coefficients $d^{\beta}+\frac{1-D}{2}$ in \EQUref{1.100.100}
are zero. 
However we could define $C^{\alpha \beta}_{\gamma}$ 
by $\Gamma^{\alpha \beta}_{\gamma}$ also for this case 
because of the property of the unit. 
Then we see that $C^{\alpha \beta}_{\gamma}$ satisfies 
the conditions of Frobenius manifold through case by case argument and 
the technique of flat pencils. 

By this argument, we have a uniqueness of the product structure 
if the unit vector field and the holomorphic metric is given 
(cf. \cite[p.195, Remark 3.2]{Dubrovin}). 
However we give a stronger uniqueness theorem 
which asserts that the Frobenius manifold structure with intersection form $I^*$ 
is unique up to $\C^*$-action 
under the condition that only the Euler vector field is fixed.

Our results have an interesting application. 
A flat coordinate system $t^1,\cdots,t^{\mu}$ 
could be constructed 
by Jacobi forms for $G_2$ case \cite{Kokyuroku} and 
for $E_6$ case \cite{E6}.  
Thus the functions $I^*(dt^i,dt^j)$ on $M$ for the symmetric tensor $I^*$ 
are represented by the derivatives of these automorphic functions. 
On the other hand, the result that $I^*$ is an intersection 
form of the Frobenius manifold means that the functions  
$I^*(dt^i,dt^j)$ have a lot of differential relations. 
Thus we could obtain many differential relations 
for these automorphic functions systematically. 
An explicit formula of differential relations 
is partly announced in \cite{D4}. 

In the results of this paper, 
the part of constructing the Frobenius manifold 
is already announced in \cite{D4} in the form 
of an existence of a potential. 
For an explicit calculation of the potential, 
it is done for $D_4$ case \cite{D4}  
and for $G_2$ case \cite{Bertola}. 

This paper is organized as follows. 

In Section 2, we review the notions necessary for 
later sections. We introduce a concept of an elliptic root system, 
the elliptic Weyl group, 
the domain $\Ee$, the symmetric tensor on $\Ee$, 
invariant ring and the Euler operator. 
Here we introduce new notions of a signed marking and 
an orientation of an elliptic root system. 
These notions fit the natural definition of the domain 
$\Ee$ which we shall study. 

In Section 3, we give the results. 
We introduce the quotient space of $\Ee$ by the elliptic Weyl group 
and see that the symmetric tensor descends to the quotient space. 
Then we assert that the quotient space has a structure of Frobenius manifold 
which is compatible with the symmetric tensor 
under some suitable condition (so-called ``codimension 1"). 
We also assert a strong uniqueness theorem of the 
structure of Fronbenius manifold. 

In Section 4, we show that the quotient space have a 
structure of Frobenius manifold. 
First we review the work of construction of holomorphic 
metric by \cite{extendedII}. 
Then we define the multiplication in a case by case manner. 
We construct the potential of this product. 
This potential is used to show the properties of the product. 
For these results, the technique of a flat pencil is 
necessary and we collect the necesary results in the beginning 
of this section. 

In Section 5, we give a strong uniqueness theorem 
which asserts that the Frobenius manifold structure with intersection form $I^*$ 
is unique up to $\C^*$-action 
under the condition that only the Euler vector field is fixed.

The auther would like to thank Prof. Michihisa Wakui 
for his careful reading of the manuscript and 
for his continuous encouragement. 

%
%
%
%
%
%
%
%
%
%
%
%

\section{Weyl group invariant ring}
The purpose of this section is to review the invariant 
ring and the Euler operator introduced in \cite{extendedII}. 
The notions of a signed marking and an orientation 
of an elliptic root system are new.
\subsection{Elliptic root system}
In this subsection, we define an elliptic root system and 
its orientation.

Let $l$ be a positive integer. 
Let $F$ be a real vector space of rank $l+2$ 
with a negative semi-definite or 
positive semi-definite 
symmetric bilinear form 
$I:F \times F \to \R$, whose radical 
$\mathrm{rad}I:=\{x \in F\,|\,I(x,y)=0,\forall y \in F\}$ 
is a vector space of rank 2. 
For a non-isotropic element $\alpha \in F$ (i.e. 
$I(\alpha,\alpha) \neq 0$), we put 
$\alpha^{\vee}:=2\alpha/I(\alpha,\alpha) \in F$. 
The reflection $w_{\alpha}$ with respect to $\alpha$ is defined by 
\begin{equation}\EQUlabel{2.11}
w_{\alpha}(u):=u-I(u,\alpha^{\vee})\alpha\quad
(\forall u \in F).
\end{equation}
%
%
\begin{Def}
{\rm 
(\cite[p.104, Def. 1]{extendedI}) }
A set $R$ of non-isotropic elements of $F$ is an elliptic 
root system belonging to $(F,I)$ if it satisfies the axioms 
1-4:
\begin{enumerate}
	\item The additive group generated by $R$ in $F$, 
	denoted by $Q(R)$, is a full sub-lattice of $F$. 
	That is, the embedding $Q(R) \subset F$ induces 
	the isomorphism : $Q(R) \otimes_{\Z}\R \simeq 
	F$. 
	\item $I(\alpha,\beta^{\vee}) \in \Z$ for 
	$\alpha,\beta \in R$.
	\item $w_{\alpha}(R)=R$ for $\forall \alpha \in R$. 
	\item If $R=R_1 \cup R_2$, with $R_1 \perp R_2$, 
	then either $R_1$ or $R_2$ is void. 
\end{enumerate}
\end{Def}
For an elliptic root system $R$ belonging to $(F,I)$, 
the additive group $\mathrm{rad}I \cap Q(R)$ 
is isomorphic to $\Z^2$. 

\begin{Def}
An elliptic root system $R$ is called oriented 
if the $\R$-vector space $\mathrm{rad}I$ is oriented. 
A frame $\{a,b\}$ of $\mathrm{rad}I$ is called 
admissible if 
$\mathrm{rad}I \cap Q(R) \simeq \Z a \oplus \Z b$ and 
it gives the orientation of $\mathrm{rad}I$. 
\end{Def}

\begin{Rem}
If an elliptic root system $(R,F,I)$ comes from vanishing 
cycles of a Milnor fiber of a simple elliptic singularity, then 
$\mathrm{rad}I \cap Q(R) \simeq H_1(E_{\infty},\Z)$ 
for an elliptic curve $E_{\infty}$ at infinity 
(cf. \cite[p.18]{extendedII}). 
Then $(R,F,I)$ is canonically oriented by the complex structure 
of the elliptic curve $E_{\infty}$. 
\end{Rem}

\subsection{Hyperbolic extension and Weyl group}
In this subsection, we define a signed marking, 
a hyperbolic extension and its Weyl group.  

\begin{Def}
Let $R$ be an elliptic root system $R$ belonging to $(F,I)$. 
By a signed marking, we mean a non-zero element $a$ of 
$\mathrm{rad}I \cap Q(R)$ such that $Q(R) \cap \R a=\Z a$.
\end{Def}
Hereafter we fix an oriented elliptic root 
system with a signed marking $(R,a)$ such that the quotient root system 
$R/\R a$ ($:=\mathrm{Image}(R \hookrightarrow 
F \to F/\R a$)) is reduced (i.e. $\alpha,c\alpha \in R/\R a$ 
implies $c \in \{\pm 1\}$). 

Let $F^1$ be a real vector space of rank $l+3$ and 
$I^1:F^1 \times F^1 \to \R$ an $\R$-symmetric bilinear form. 
The pair $(F^1,I^1)$ is called a hyperbolic extension of 
$(F,I)$ if 
$F^1$ contains $F$ as a linear subspace, 
$\mathrm{rad}I^1=\R a$ and $I^1|_{F}=I$. 
A hyperbolic extension is unique up to isomorphism. 
Hereafter we fix a hyperbolic extension $(F^1,I^1)$.

We define a reflection 
$\tilde{w}_{\alpha} \in GL(F^1)$ 
by $\tilde{w}_{\alpha}(u):=u-I^1(u,\alpha^{\vee})\alpha$ 
for $u \in F^1$.

We define a Weyl group $\tilde W$ 
(resp. $W$) by 
\begin{equation}
\tilde W:=\langle \tilde{w}_{\alpha}\,|\,\alpha \in R\rangle\ 
(\mbox{resp. }W:=\langle{w}_{\alpha}\,|\,\alpha \in R\rangle).
\end{equation}

We have a natural exact sequence:
\begin{equation}
0 \to K_{\Z} \to \tilde{W} \to W \to 1,
\end{equation}
where $\tilde{W} \to W$ is given by 
the restriction of $\tilde{W}$ on $F$ and 
$K_{\Z}$ is the kernel of $\tilde{W} \to W$. 
The group $K_{\Z}$ is isomorphic to $\Z$. 

\subsection{Domain}
In this subsection, we define a domain $\Ee$ 
for the oriented elliptic root system with 
the signed marking $(R,a)$ belonging to $(F,I)$ 
such that $R/\R a$ is reduced. 

For $(F,I)$, the set 
$\{c \in \R\,|\,$ the bilinear form $cI$ defines a semi-positive 
even lattice structure on $Q(R)\}$ has the unique element of 
the smallest absolute value. 
We denote it by $(I_R:I)$.

Take $b \in \mathrm{rad}I \cap Q(R)$ 
such that $\{a,b\}$ gives an admissible frame. 
Then we could choose 
an isomorphism $\varphi:\Z \simeq K_{\Z}$ 
and $\tilde{\lambda} \in F^1 \setminus F$ 
such that 
\begin{equation}\EQUlabel{2.3.1}
\varphi(n) \cdot \tilde{\lambda}=\tilde{\lambda}-na\ 
(n \in \Z),\quad
(I_R:I)I(\tilde{\lambda},b)>0.
\end{equation}
By the condition (\ref{2.3.1}), 
$\tilde{\lambda}$ is unique up to adding an element of $F$, 
and such an isomorphism $\varphi$ is unique. 

We define two domains:
\begin{align}
	\Ee&:=\{
	x \in \mathrm{Hom}_{\R}(F^1,\C)\,|\,
	\langle a,x \rangle=1,\ 
	\mathrm{Im}\langle b,x \rangle>0\ \},\\
	\Hh&:=\{
	x \in \mathrm{Hom}_{\R}(\mathrm{rad}I,\C)\,|\,
	\langle a,x \rangle=1,\ 
	\mathrm{Im}\langle b,x \rangle>0\ \},
\end{align}
where $\langle\ ,\  \rangle$ is the natural pairing 
$F^1_{\C}\times (F^1_{\C})^* \to \C$ and $F^1_{\C}:=F^1\otimes_{\R}{\C}$. 
We have a natural projection 
\begin{equation}
\pi:\Ee \to \Hh.
\end{equation}
\indent
For a root $\alpha \in R$, we define the reflection hyperplane of 
$\Ee$ by 
\begin{equation}
H_{\alpha}:=\{x \in \Ee\,|\,\langle \alpha,x \rangle=0\,\}.
\end{equation}

We define a left action of $\tilde{W}$ on $\Ee$ 
by 
\begin{equation}
\langle \eta,g \cdot x \rangle
:=
\langle g^{-1}\eta,x \rangle
\end{equation}
for $g \in \tilde{W},\ \eta \in F^1_{\C},\ x \in \Ee$. 

For a complex manifold $M$, 
we denote by $\Oo_M$ 
(resp. $\Omega^1_M,\ \Theta_M$)
the sheaf of holomorphic functions 
(resp. holomorphic 1-forms, holomorphic vector fields). 
We denote 
by $\Oo(M)$ (resp. $\Omega^1(M),\ \Theta(M)$) 
the module $\Gamma(M,\Oo_M)$ 
(resp. $\Gamma(M,\Omega^1_M),\ \Gamma(M,\Theta_M)$).

We define a vector field $E'$ on $\Ee$ by the conditions 
$$
E'x=0\quad (\forall x \in F),\quad
E'\tilde{\lambda}=\frac{1}{-2\pi\sqrt{-1}}.
$$
The vector field $E'$ is uniquely determined by the condition (\ref{2.3.1}). 

We define a $\C$-symmetric bilinear form $I^*$ on $\Omega^1_{\Ee}$. 
Since we have a canonical isomorphism 
$T^*_p\Ee \simeq \C \otimes_{\R}(F^1/\R a)$ for $p \in \Ee$, 
we have an $\Oo_{\Ee}$-bilinear form
\begin{equation}
I^*:\Omega^1_{\Ee}\times \Omega^1_{\Ee}\to \Oo_{\Ee}
\end{equation}
induced from 
$I^1:F^1/\R a\times F^1/\R a \to \R$. 
We remark that $Lie_{E'} I^*=0$, where $Lie$ is the Lie derivative.

\newcommand{\EEE}{(\Ee \setminus \cup_{\alpha \in R}H_{\alpha})/\tilde{W}}
The action of $\tilde W$ on $\Ee$ 
(resp. $\Ee \setminus \cup_{\alpha \in R}H_{\alpha}$) 
is properly discontinuous 
(resp. properly discontinuous and fixed point free) (cf.\cite{extendedII}), 
thus the quotient space $\Ee/\tilde{W}$ 
(resp. $\EEE$) 
has a structure of analytic space 
(resp. complex manifold). 
We have the following diagram of analytic spaces:
\begin{equation}
\begin{CD}
\EEE @>i_1>> \Ee/\tilde{W} \\
@V{\pi}VV @V{\pi}VV \\
\Hh @= \Hh 
\end{CD}.
\end{equation}
The morphism $i_1$ is an open immersion and its image is open dense. 

Since the tensors $E'$ and $I^*$ on $\Ee$ are 
$\tilde{W}$-invariant, 
these tensors descend to the space 
$\EEE$:
\begin{align}
E':\ &\Omega^1_{\EEE} \to \Oo_{\EEE},\EQUlabel{2.3.1001}\\
I^*:\ &\Omega^1_{\EEE}\times \Omega^1_{\EEE} \to \Oo_{\EEE}.
\EQUlabel{2.3.1002}
\end{align}
\subsection{The Weyl group invariant ring}
In order to extend the domain of the definition of the tensors 
\EQUref{2.3.1001} and \EQUref{2.3.1002}, 
we introduce the Weyl group invariant ring $\Ss^W$, 
$\Ss^W$ modules $\Omega^1_{\Ss^W}, \Der_{\Ss^W}$ 
in this subsection and formulate the 
tensors $E'$ and $I^*$ by these $\Ss^W$ modules 
in the next subsection. 

We define an $\Oo_{\Hh}$-module $\Ss^W_k$ of 
$\tilde{W}$-invariant functions parametrized by $k \in \Z$ 
as the subsheaf of $\pi_*\Oo_{\Ee}$ by 
\begin{equation}
\Ss^W_k(U):=\{
f\in  \pi_*\Oo_{\Ee}(U)\,|\,
f(g \cdot x)=f(x) \ (\forall g \in \tilde{W}, \forall x \in 
\pi^{-1}(U)),\ 
E'f=kf
\}
\end{equation}
for an open set $U \subset \Hh$.

We define the $\Oo_{\Hh}$-graded algebra $\Ss^W$ by 
\begin{equation}
\Ss^W:=\bigoplus_{k \in \Z}\Ss^W_k.
\end{equation}
We have injective homomorphisms
\begin{equation}
\Ss^W \to \pi_*\Oo_{\Ee/\tilde{W}} \to \pi_* \Oo_{\EEE}.
\end{equation}

\begin{Thm}\EQUlabel{2.10}
{\rm(
\cite{Chevalley3}, 
\cite{Chevalley4}, 
\cite{Chevalley6},
\cite{Chevalley2}, 
\cite{Chevalley1}, 
\cite{Chevalley5}
)}

The $\Oo_{\Hh}$-graded algebra $\Ss^W$ is an 
$\Oo_{\Hh}$-free algebra, i.e. 
\begin{equation}
\Ss^W=\Oo_{\Hh}[s^1,\cdots,s^{n-1}]
\end{equation}
for $s^j \in \Ss^W_{c^j}(\Hh)$ with 
$c^1 \geq c^2 \geq \cdots \geq c^{n-1}>0$ 
and $n:=l+2$. 
We remark that $j$ of $s^j$ is a suffix.
\end{Thm}
We introduce the $n$-th invariant $s^n \in \Ss^W_{0}(\Hh)$. 
Since an element of $\mathrm{rad}I$ naturally gives an element of 
$\Ss^W_0(\Hh)$, we define $s^n=b$, 
where $b \in \mathrm{rad}I$ is introduced in Section 2.3. 
We put $c^n=0$. Then $s^j \in \Ss^W_{c^j}(\Hh)$ for $j=1,\cdots,n$. 

We define two $\Ss^W$-modules $\Omega^1_{\Ss^W}$ 
and $\Der_{\Ss^W}$. 

For an open set $U \subset \Hh$, we put
\begin{equation}
\Omega^1_{\Ss^W}(U):=
\Omega^1_{\Ss^W(U)/\C},
\end{equation}
where R.H.S. is the module of relative differential forms 
of a $\C$-algebra $\Ss^W(U)$. 
Since $\Ss^W$ is an $\Oo_{\Hh}$-free algebra, 
$\Omega^1_{\Ss^W}$ defines a sheaf. 
$\Omega^1_{\Ss^W}$ has an $\Ss^W$-module structure. 
We remark that since a local section of $\Omega^1_{\Ss^W}$ 
determines a local section $\pi_*\Omega^1_{\Ee}$ and 
$\pi_*\Omega^1_{\Ee \setminus \cup_{\alpha \in R}H_{\alpha}}$ 
which is $\tilde{W}$-invariant, 
there exists a natural lifting map 
$\Omega^1_{\Ss^W} \to \pi_*\Omega^1_{\EEE}$. 

We put 
\begin{equation}
\Der_{\Ss^W}:=
\underline{\mathrm{Hom}}_{\Ss^W}(\Omega^1_{\Ss^W},\Ss^W).
\end{equation}

By the generators of $\Ss^W$, we have
\begin{equation}\EQUlabel{2.4.3000}
\Omega^1_{\Ss^W}=\oplus_{i=1}^n \Ss^W ds^i,\quad
\Der_{\Ss^W}=\oplus_{i=1}^n \Ss^W \frac{\p}{\p s^i}.
\end{equation}

\subsection{Euler operator and bilinear form on the invariant ring}
We define the Euler operator $E$ 
as a vector field on $\Ee$ defined by 

\begin{equation}
E:=\frac{1}{c^1}E'.\EQUlabel{2.5.999}
\end{equation}
As in the case of $E'$ in \EQUref{2.3.1001}, 
$E$ defines a morphism
\begin{equation}
E:\ \Omega^1_{\EEE} \to \Oo_{\EEE}.
\EQUlabel{2.5.998}
\end{equation}

By \cite{extendedII}, we have the 
$\Ss^W$-homomorphism and $\Ss^W$-symmetric bilinear form:
\begin{align}
E:\ &\Omega^1_{\Ss^W} \to \Ss^W\EQUlabel{2.5.1003},\\
I^*:\ &\Omega^1_{\Ss^W}\times \Omega^1_{\Ss^W} \to \Ss^W
\EQUlabel{2.5.1004}
\end{align}
with the following diagrams:
\begin{equation}
\begin{CD}
\pi_*\Omega^1_{\EEE} @>E>> \pi_*\Oo_{\EEE}\\
@AAA @AAA\\
\Omega^1_{\Ss^W} @>E>> \Ss^W
\end{CD},\EQUlabel{2.5.1001}
\end{equation}
\begin{equation}
\begin{CD}
	I^*:\  & 
	\pi_*\Omega^1_{\EEE} & \times & 
	\pi_*\Omega^1_{\EEE} & @>>> 
	\pi_*\Oo_{\EEE}\\
	& @AAA @AAA & @AAA\\
	I^* :\  & \Omega^1_{\Ss^W} & \times & \Omega^1_{\Ss^W}&
	@>>> \Ss^W,
\end{CD}\EQUlabel{2.5.1002}
\end{equation}
where the upper line of \EQUref{2.5.1001} is induced by \EQUref{2.5.998} 
and the upper line of \EQUref{2.5.1002} is induced by \EQUref{2.3.1002}. 
The morphisms 
\EQUref{2.5.1003} and \EQUref{2.5.1004}
are uniquely characterized by the diagrams 
\EQUref{2.5.1001} and \EQUref{2.5.1002}
respectively because 
$
\Ss^W \to \pi_*\Oo_{\EEE}
$
is injective.

\section{Results}
In this section, we first define the Weyl group quotient space $\EW$
for $\Ee$ and $\tilde{W}$ defined in Section 2. 
Then we assert that $\EW$ has a structure of Frobenius manifold 
under some suitable condition. 

\subsection{Weyl group quotient space}
In this subsection, we define the Weyl group quotient space 
and study the tensors on it. 

Let $\Ee$, $\Hh$ be the domains and $\tilde{W}$ be
the Weyl group defined in Section 2. 
Let $(An)$ and $(Set)$ be categories of analytic spaces and 
sets, respectively. Let $((An)/\Hh)^{\circ}$ be 
the dual category of the category of $\Hh$-objects. 
Since the $\Oo_{\Hh}$-algebra $\Ss^W$ is of finite presentation 
(Theorem \ref{2.10}), 
the analytic space $\Specan\,\Ss^W$ could be defined by \cite{Houzel}. 
We define the Weyl group quotient space $\EW$ by 
\begin{equation}
\EW:=\Specan\,\Ss^W.\EQUlabel{3.1.20}
\end{equation}
We denote the structure morphism 
$\EW \to \Hh$ also by $\pi$. 
The space $\EW$ is isomorphic to $\Hh \times \C^{n-1}$ by 
Theorem \ref{2.10}.

By definition of $\Specan$, there exists a 
natural isomorphism:
\begin{equation*}
\mathrm{Hom}_{(An)/\Hh}(X,\EW)\simeq 
\mathrm{Hom}_{\Oo_{X}}(f^*\Ss^W,\Oo_X)
\end{equation*}
for an object $f:X \to \Hh$ of the category $(An)/\Hh$. 
Since there exists 
a canonical isomorphism: 
$
\mathrm{Hom}_{\Oo_{X}}(f^*\Ss^W,\Oo_X)
\simeq 
\mathrm{Hom}_{\Oo_{\Hh}}(\Ss^W,f_*\Oo_X)
$, 
we have
\begin{equation}
\mathrm{Hom}_{(An)/\Hh}(X,\EW)
\simeq 
\mathrm{Hom}_{\Oo_{\Hh}}(\Ss^W,f_*\Oo_X).
\EQUlabel{3.1.10}
\end{equation}

We define a ringed space $(\Hh,\Ss^W)$ 
by the space $\Hh$ with the sheaf $\Ss^W$. 
We define a morphism of the category of ringed spaces:
\begin{equation}
\varphi:(\EW,\Oo_{\EW}) \to (\Hh,\Ss^W)\EQUlabel{3.1.3001}
\end{equation}
by the mapping $\pi:\EW \to \Hh$ and the morphism 
\begin{equation}
\phi:\Ss^W \to \pi_*\Oo_{\EW} \EQUlabel{3.1.2000}
\end{equation}
which corresponds to the identity element of 
$\mathrm{Hom}_{(An)/\Hh}(\EW,\EW)$ by \EQUref{3.1.10}. 

\begin{Prop}
We have the canonical isomorphism:
\begin{equation}\EQUlabel{3.1.4001}
\varphi^*\Omega^1_{\Ss^W} \simeq \Omega^1_{\EW}.
\end{equation}
\end{Prop}
\begin{proof}
We define the ringed space $(\EW^{alg},\Oo_{\EW^{alg}})$ as follows: 
As a set, $\EW^{alg}=\EW$. 
A topology on $\EW^{alg}$ is introduced so that 
$\{(U,f)\subset \EW^{alg}\,|\,U \subset \Hh:\hbox{open},
f \in \Gamma(U,\Ss^W)\}$ 
becomes an open basis, 
where 
$(U,f):=\{x \in \EW^{alg}\,|\,\pi(x) \in U,f(x) \neq 0\,\}$. 
We define the sheaf $\Oo_{\EW^{alg}}$ associated with the presheaf 
$\Oo_{\EW^{alg}}((U,f)):=\Gamma(U,\Ss^W)_f$ 
for an open set $(U,f)$. 

The morphism $\varphi:(\EW,\Oo_{\EW}) \to (\Hh,\Ss^W)$ 
factors as the composite of the morphisms:
\begin{equation}
(\EW,\Oo_{\EW})
\stackrel{\varphi_1}{\to}
(\EW^{alg},\Oo_{\EW^{alg}})
\stackrel{\varphi_2}{\to}
(\Hh,\Ss^W).
\end{equation}
We define a sheaf $\Omega^1_{\EW^{alg}}$ on $\EW^{alg}$ 
as a sheafification of the presheaf
\begin{equation}
(U,f) \mapsto \Omega^1_{\Ss^W}(U)_f.
\end{equation}
Then we have a natural isomorphism 
\begin{equation}
\varphi_2^*\Omega^1_{\Ss^W}\simeq \Omega^1_{\EW^{alg}}
\end{equation}
by a discussion of an affine morphism in scheme theory. 
Also we have a natural isomorphism 
\begin{equation}
\varphi_1^*\Omega^1_{\EW^{alg}}\simeq
\Omega^1_{\EW}
\end{equation}
because an $\Oo_{\EW^{alg}}$-locally free basis 
of algebraic 1-forms is regarded as 
an $\Oo_{\EW}$-locally free basis 
of analytic 1-forms. 
\end{proof}

We define the $\Oo_{\EW}$-homomorphism and 
$\Oo_{\EW}$-symmetric bilinear form 
\begin{align}
E:\ &\Omega^1_{\EW} \to \Oo_{\EW},\EQUlabel{3.1.1001}\\
I^*_{\EW}:\ &\Omega^1_{\EW} \times \Omega^1_{\EW} \to \Oo_{\EW},
\EQUlabel{3.1.1002}
\end{align}
by taking the pull-back of \EQUref{2.5.1003} and \EQUref{2.5.1004}
by $\varphi$. 

We shall see the relation between 
\EQUref{3.1.1001} (resp.\EQUref{3.1.1002}) on $\EW$ and 
\EQUref{2.5.998} (resp.\EQUref{2.3.1002}) on $\EEE$.

By \EQUref{3.1.10}, a natural inclusion 
$\Ss^W \hookrightarrow \pi_*\Oo_{\Ee/\tilde{W}}$ corresponds 
to the mapping 
\begin{equation}
i_2:\Ee/\tilde{W} \to \EW.
\end{equation}
By \cite{extendedII}, the morphism $i_2$ is an open immersion and we have 
$\EW \simeq \Ee/\tilde{W}\cup \Hh$. 
The composite mapping $\EEE \stackrel{i_1}{\to} \Ee/\tilde{W}
\stackrel{i_2}{\to}\EW$ is also an open immersion and 
its image is open dense. 

We have the following diagram of ringed spaces:
\begin{equation}
\begin{CD}
(\EEE,\Oo_{\EEE}) @>i_1>>
(\Ee/\tilde{W},\Oo_{\Ee/\tilde{W}}) @>i_2>>
(\EW,\Oo_{\EW})\\
@V{\pi}VV @V{\pi}VV @V{\varphi}VV\\
(\Hh,\Oo_{\Hh}) @=
(\Hh,\Oo_{\Hh}) @<<<
(\Hh,\Ss^W)
\end{CD}.
\end{equation}
\begin{Prop}\EQUlabel{3.1.1500}
The $\Oo_{\EEE}$-homomorphism \EQUref{2.5.998} 
(resp. the $\Oo_{\EEE}$-symmetric bilinear form \EQUref{2.3.1002})
is uniquely extended to 
the $\Oo_{\EW}$-homomorphism 
$\Omega^1_{\EW} \to \Oo_{\EW}$ 
(resp. $\Oo_{\EW}$-symmetric bilinear form 
$\Omega^1_{\EW}\times \Omega^1_{\EW} \to \Oo_{\EW}$) 
and coincides with 
\EQUref{3.1.1001}
(resp.\EQUref{3.1.1002}). 
\end{Prop}
\begin{proof}
Since the image of the open immersion $i_2 \circ i_1$ 
is open dense, 
we should only prove that 
the pull-back of \EQUref{3.1.1001}(resp.\EQUref{3.1.1002})
by $i_2 \circ i_1$ coincides 
with \EQUref{2.5.998}(resp.\EQUref{2.3.1002}). 
The former is the pull-back of 
\EQUref{2.5.1003}(resp.\EQUref{2.5.1004}) 
by $\varphi \circ i_2 \circ i_1$. 
The latter could be written as 
$E:(\varphi \circ i_2 \circ i_1)^*\Omega^1_{\Ss^W} \to 
(\varphi \circ i_2 \circ i_1)^*\Ss^W$ 
(resp. 
$I^*:
(\varphi \circ i_2 \circ i_1)^*\Omega^1_{\Ss^W}
\times
(\varphi \circ i_2 \circ i_1)^*\Omega^1_{\Ss^W}
\to
(\varphi \circ i_2 \circ i_1)^*\Ss^W$). 

Then we have the result by applying the 
following lemma to 
\EQUref{2.5.1001} and \EQUref{2.5.1002} 
using the fact that $(\varphi \circ i_2 \circ i_1)_*=\pi_*$. 
\end{proof}

\begin{Lemma}
Let $f:(X,\Oo_X)\to (Y,\Oo_Y)$ be a morphism 
of ringed spaces. 
Let $\mathcal{F},\mathcal{G}$ be $\Oo_Y$-modules. 
If we have $\alpha:f^*\mathcal{F} \to f^*\mathcal{G}$, 
$\beta:\mathcal{F}\to \mathcal{G}$ and 
a commutative diagram:
\begin{equation}
\begin{CD}
f_*f^*\mathcal{F} @>{f_*\alpha}>> f_*f^*\mathcal{G}\\
@AAA @AAA\\
\mathcal{F} @>{\beta}>> \mathcal{G}
\end{CD}
\end{equation}
for the natural morphisms $\mathcal{F}\to f_*f^*\mathcal{F}$ 
and $\mathcal{G}\to f_*f^*\mathcal{G}$, 
then we have $\alpha=f^*\beta$. 
\end{Lemma}
\begin{proof}
By the naturality of $f^*f_* \to id.$, we have the commutative diagram:
\begin{equation}
\begin{CD}
f^*\mathcal{F} @>{\alpha}>> f^*\mathcal{G}\\
@AAA @AAA\\
f^*f_*f^*\mathcal{F} @>{f^*f_*\alpha}>> f^*f_*f^*\mathcal{G}\\
@AAA @AAA\\
f^*\mathcal{F} @>{f^*\beta}>> f^*\mathcal{G}. 
\end{CD}
\end{equation}
Since the composite morphism 
$f^*\mathcal{F} \to f^*f_*f^*\mathcal{F} \to f^*\mathcal{F}$ 
is the identity morphism, we have the result. 
\end{proof}

\subsection{Frobenius manifold}
In this section, we give the main theorem which asserts 
that the space $\EW$ admits a structure of Frobenius manifold 
and it is unique up to $\C^*$ action under some suitable condition. 

We first remind the definition of Frobenius manifold and its intersection 
form.
\begin{Def}\rm{(\cite[p.146, Def. 9.1]{Hertling})}\EQUlabel{3.333}
A Frobenius manifold is a tuple $(M,\circ,e,E,g)$ 
where $M$ is a complex manifold of dimension $\geq 1$ 
with holomorphic metric $g$ and multiplication $\circ$ 
on the tangent bundle, $e$ is a global unit field 
and $E$ is another global vector field, subject to the following 
conditions:
\begin{enumerate}
\item the metric is invariant under the multiplication, i.e., 
$g(X\circ Y,Z)=g(X,Y\circ Z)$
for local sections $X,Y,Z \in \Theta_M$,
\item (potentiality) the $(3,1)$-tensor 
$\nabla \circ $ is symmetric
(here, $\nabla$ is the Levi-Civita connection of the metric), i.e., 
$
\nabla_X(Y \circ Z)-Y \circ \nabla_X(Z)
-\nabla_Y(X \circ Z)+X \circ \nabla_Y(Z)-[X,Y]\circ Z=0,
$ for local sections $X,Y,Z \in \Theta_M$,
\item the metric $g$ is flat,
\item $e$ is a unit field and it is flat, i.e. $\nabla e=0$, 
\item the Euler field $E$ satisfies 
$Lie_{E}(\circ )=1 \cdot \circ$ and 
$Lie_{E}(g)=D \cdot g$ for some $D \in \C$.
\end{enumerate}
\end{Def}
\begin{Def}\rm{(\cite[p.191]{Dubrovin})}\EQUlabel{3.555}
For a Frobenius manifold $(M,\circ,e,E,g)$, 
we define an intersection form 
$h^*:\Omega^1_M \times \Omega^1_M \to \Oo_M$ 
by 
\begin{equation}\EQUlabel{2.9.3}
h^*(\omega_1,\omega_2)=
	g(E,g^*(\omega_1)\circ g^*(\omega_2))
\end{equation}
where $g^*:\Omega^1_{M} \to \Theta_{M}$
is the isomorphism induced by $g$. 
\end{Def}

For the oriented elliptic root system 
with the signed marking $(R,a)$ such that $R/\R a$ is reduced, 
the condition $c^1 >c^2$ is called ``codimension 1" in \cite{extendedII}. 
\begin{Thm}\EQUlabel{2.9.1}
If the oriented elliptic root system with the signed marking $(R,a)$
such that $R/\R a$ is reduced satisfies the condition of codimension 1, 
then we have the following results.\\
(1) $\EW$ has a structure of 
Frobenius manifold $(\EW,\circ,e,E,J)$ with the following conditions:
\begin{enumerate}
	\item $E$ is the Euler field defined in \EQUref{3.1.1001}.
	\item $I^*_{\EW}$ gives the intersection form 
	of a Frobenius manifold $(\EW,\circ,e,E,J)$. 
\end{enumerate}
(2) For $c \in \C^*$, $(\EW,c^{-1}\circ,ce,E,c^{-1}J)$ is also a 
Frobenius manifold satisfying conditions of (1). \\
(3) Let $(\EW,\circ',e',E',J')$ be a Frobenius manifold 
which satisfies conditions of (1). Then there exists $c \in \C^*$ 
such that 
$(\EW,\circ',e',E',J')=(\EW,c^{-1}\circ,ce,E,c^{-1}J)$. 
\end{Thm}
\begin{Rem}
In the definition of Frobenius manifold, the homogeneity 
$Lie_E(J)=D \cdot J$ for some $D \in \C$ is assumed. 
By the equation \EQUref{2.9.3}, $D$ must be $1$ because 
$Lie_E I^*_{\EW}=0$ and $Lie_E(\circ)=\circ$. 
\end{Rem}
\begin{Rem}
The existence of the holomorphic metric $J$ is 
already shown in \cite{extendedII}. 
The existence of the product sturcture is 
already announced in \cite{D4} 
in the form of the existence of the potential. 
\end{Rem}

\section{Construction of Frobenius manifold structure}
In this section, we give the proof of Theorem \ref{2.9.1}(1), 
that is,
the existence of Frobenius manifold structure on $\EW$.
In Section 4.1, we review a construction \cite{extendedII} 
of a flat metric on $\EW$ (Proposition \ref{4.16.1}) and flat coordinates. 
In Section 4.2, we recall the notion of a flat pencil. 
In Section 4.3, we construct a product structure 
on the tangent space of $\EW$. 
In Section 4.4, we construct a potential of 
the product. 
In Section 4.5, we show the properties 
of the product. 
In Section 4.6, we show that these constructions 
give a Frobenius manifold structure. 

Hereafter we shall calculate tensors by using indices. 
In that case, we use Einstein's summation convention, that is, 
if an upper index of one tensor and a lower of the other tensor 
coincide, then we take summation for the same letter.

\subsection{A construction of a flat metric and flat coordinates}
Let $\EW$ be the Weyl group quotient space defined in \EQUref{3.1.20}.
Hereafter we assume that $(R,a)$ is codimension 1.

We prepare the relation between $\Ss^W$-modules and $\Oo_{\EW}$-modules. 

We first define a notion of degree. 
For $f \in \Ss^W$ and $d \in \Q$, if $Ef=df$, then we call $d$ 
{\it the degree of $f$}. 
For $f \in \Ss^W_k$, the degree of $f$ is $\frac{k}{c^1}$. 
Especially the degree of $s^i$ in Theorem \ref{2.10} is 
$d^i:=\frac{c^i}{c^1}\ (i=1,\cdots,n)$, i.e. 
\begin{equation}
E s^i=d^i s^i\ (i=1,\cdots,n), \quad
1=d^1>d^2 \geq \cdots \geq d^{n-1} > d^n=0.
\EQUlabel{2.5.1}
\end{equation}
A degree is defined also for local sections of 
$\Omega^1_{\Ss^W}$ and $\Der_{\Ss^W}$. They have $\Ss^W$-free 
homogeneous generators by \EQUref{2.4.3000}.

We have morphisms:
\begin{align}
\Ss^W&\to \varphi_*\varphi^*\Ss^W
\simeq \varphi_*\Oo_{\EW},\EQUlabel{4.1.3001}\\
\Omega^1_{\Ss^W}&\to \varphi_*\varphi^*\Omega^1_{\Ss^W}
\simeq \varphi_*\Omega^1_{\EW},\EQUlabel{4.1.3002}\\
\Der_{\Ss^W}&\to \varphi_*\varphi^*\Der_{\Ss^W}
\simeq \varphi_*\Theta_{\EW}.\EQUlabel{4.1.3003}
\end{align}
where \EQUref{4.1.3001} is $\phi:\Ss^W \to \pi_*\Oo_{\EW}$ 
defined in \EQUref{3.1.2000}. 
The morphisms \EQUref{4.1.3002} and \EQUref{4.1.3003} are defined for a 
morphism $\varphi:(\EW,\Oo_{\EW}) \to (\Hh,\Ss^W)$.  

\begin{Prop}\EQUlabel{3.1.5}
The morphisms \EQUref{4.1.3001}, \EQUref{4.1.3002} and \EQUref{4.1.3003} are 
injective. A homogeneous local section 
with respect to the Euler operator $E \in \Theta(\EW)$ of 
$\varphi_*\Oo_{\EW}$ 
(resp. $\varphi_*\Omega^1_{\EW}$, $\varphi_*\Theta_{\EW}$) 
is an image of 
$\Ss^W$ (resp. $\Omega^1_{\Ss^W}$, $\Der_{\Ss^W}$).
\end{Prop}
\begin{proof}
The morphism $\varphi$ decomposes into $\varphi=\varphi_2 \circ \varphi_1$ 
as in the proof of proposition \ref{3.1.4001}. 
For a $\Ss^W$-module $\mathcal{M}$, 
we have $\mathcal{M}\simeq {\varphi_2}_*\varphi_2^*\mathcal{M}$ 
because $\varphi_2$ is an analogue of affine morphism 
of scheme theory. 
Also $\varphi_2^*\mathcal{M} \to {\varphi_1}_*\varphi_1^*
(\varphi_2^*\mathcal{M})$ is injective because 
$\varphi_1$ is faithfully flat by \cite{GAGA}. 
Thus we obtain the injectivity of 
$\mathcal{M} \to \varphi_*\varphi^*\mathcal{M}$. 

By semi-positivity of the degrees of $s^1,\cdots,s^n$, 
a homogeneous section of $\varphi_*\Oo_{\EW}$ is 
an image of \EQUref{4.1.3001}. 
Since $\Omega^1_{\Ss^W}$ (resp. 
$\Der_{\Ss^W}$) is a $\Ss^W$-free with homogeneous generator
$ds^1,\cdots,ds^n$ (resp. $\frac{\p}{\p s^1},\cdots,\frac{\p}{\p s^n}$), 
the morphism \EQUref{4.1.3002} (resp. \EQUref{4.1.3003}) is 
written as 
\begin{align}
\oplus_{i=1}^n \Ss^W ds^i &\to 
\oplus_{i=1}^n \varphi_*\Oo_{\EW}ds^i,\\
\hbox{(resp.} \oplus_{i=1}^n \Ss^W \frac{\p}{\p s^i} &\to 
\oplus_{i=1}^n \varphi_*\Oo_{\EW}\frac{\p}{\p s^i}).
\end{align}
Then the assersion is obvious.
\end{proof}

We prepare the notations. 
We put $S^W_k:=\Ss^W_k(\Hh),\ S^W:=\Ss^W(\Hh)$. 
Then $S^W$ is an $\Oo(\Hh)$-free algebra:
\begin{equation}
S^W=\Oo({\Hh})[s^1,\cdots,s^{n-1}].
\end{equation}
We define $S^W$-modules:
\begin{align}
\Der_{S^W}&:=\Der_{\Ss^W}(\Hh),\\
\Omega^1_{S^W}&:=\Omega^1_{\Ss^W}(\Hh).
\end{align}

We put 
\begin{align}
&\Der_{S^W}^{lowest}:=\{\de \in \Theta(\EW)\,|\,[E,\de]=-\de,
\ \de\mbox{ is non-singular }\},\\
&\Omega_{\de}:=
\{\omega \in \Omega^1_{\EW}\,|\,Lie_{\de}\omega=0\,\}\ 
\mbox{ for }\de \in \Der_{S^W}^{lowest},\EQUlabel{4.501}\\
&V:=\{\de \in \Der_{S^W}^{lowest}\,|\,\de^2I^*_{\EW}(\omega,\omega')=0\ 
,\ \forall \omega,\omega' \in \Omega_{\de}\,\EQUlabel{4.500}\}.
\end{align}

Using generators $s^1,\cdots,s^n$ in Theorem\,\ref{2.10}, 
we have 
$
\Der_{S^W}^{lowest}=\Oo^*(\Hh)\frac{\p}{\p s^1}
$ by proposition \ref{3.1.5}. 
In \cite{extendedII} it is shown that 
$V$ is non-empty and 
for any $\de \in V$, 
\begin{equation}
V=\C^* \de.\EQUlabel{4.2.1}
\end{equation}

The following proposition gives a flat metric on $\EW$.
\begin{Prop}\EQUlabel{4.16.1}
\rm{(\cite{extendedII})}
Take an arbitrary element $\widehat{e}$ of $V$. 
Then there exists a unique non-degenerate 
$\Ss^W$-symmetric bilinear form 
\begin{equation}
\widehat{J}:\Der_{\Ss^W} \times \Der_{\Ss^W} \to \Ss^W\EQUlabel{4.16.4}
\end{equation}
and its pull-back of \EQUref{4.16.4} by $\varphi$ in 
\EQUref{3.1.3001}:
\begin{equation}
\widehat{J}:\Theta_{\EW} \times \Theta_{\EW} \to \Oo_{\EW}.
\EQUlabel{4.16.2005}
\end{equation}
which is a non-degenerate 
$\Oo_{\EW}$-symmetric bilinear form. 

They are characterized by the property 
\begin{equation}\EQUlabel{4.16.2}
\widehat{J}^*(\omega_1,\omega_2)=
\widehat{e}I^*_{\EW}(\omega_1,\omega_2)
\end{equation}
for the dual metric of \EQUref{4.16.2005} and 
$\omega_1,\omega_2 \in \Omega^1_{\widehat{e}}$. 

$\widehat{J}$ is homogeneous of degree $1$, i.e. 
$Lie_{E}(\widehat{J})=\widehat{J}$. 
Furthermore, the Levi-Civita connection $\nabla^{\widehat{J}}$ 
for $\widehat{J}$ is flat and $\nabla^{\widehat{J}} \widehat{e}=0$. 
\end{Prop}

We introduce flat coordinates. Since $\EW$ is simply-connected, 
we could take functions whose differential are flat with respect to 
$\widehat{J}$. In Lemma \ref{4.150}, we show that 
they generate the ring $\Ss^W$, thus they give global coordinates for $\EW$. 

\begin{Lemma}\EQUlabel{4.150}
(1) There exist holomorphic functions 
$t^1,\cdots,t^n \in S^W$ such that \\
	(i) $\{dt^1,\cdots,dt^n\}$ gives a $\C$-basis of 
	flat sections of $\Omega^{1}_{\EW}$ with respect to $\widehat{J}$ 
	on $\EW$. \\
	(ii) $t^1,\cdots,t^n$ are homogeneous elements 
	of $S^W$ with degree $d^i$ (i.e. $Et^i=d^i t^i$), 
	where $d^i$ is defined in \EQUref{2.5.1}.\\
	(iii) $t^n=s^n$, where $s^n$ is defined after Theorem \ref{2.10}.\\
	(iv) $\widehat{e}=\frac{\p}{\p t_1}$. \\
(2) For $t^1,\cdots,t^n$, we have the following results:\\
	(i) $\Ss^W=\Oo_{\Hh}[t^1,\cdots,t^{n-1}]$. \\
	(ii) $t^1,\cdots,t^n$ give global coordinates on $\EW$. \\
	(iii) $\Omega^1_{S^W}=\oplus_{\alpha=1}^n S^W dt^{\alpha}$. 
	We remark that we use Greek letter for the suffix. \\
	(iv) We prepare elements $\frac{\p}{\p t^{\alpha}} \in \Der_{S^W}$ 
	by $\frac{\p}{\p t^{\alpha}}t^{\beta}=\de^{\beta}_{\alpha}$. 
	Then we have $\Der_{S^W}
	=\oplus_{\alpha=1}^n S^W \frac{\p}{\p t_{\alpha}}$.\\
\end{Lemma}
\begin{proof}
Since $\EW$ is simply-connected, the space 
\begin{equation}
H_1:=\{\omega \in \Gamma(\EW,\Omega^1_{\EW})\,|\,
\nabla^{\widehat{J}}\omega=0\}. 
\end{equation}
is $n$-dimensional. 
We see that any element of $H_1$ is closed 
because $\nabla^{\widehat{J}}$ is torsion-free. 
Since $Lie_E \widehat{J}=\widehat{J}$, 
a tensor $\nabla^{\widehat{J}}E$ is flat (\cite[p.147]{Hertling}). 
Then $E$ acts on $H_1$. Thus $H_1$ is identified with 
\begin{equation}
H_2:=\{\omega \in \Gamma(\Hh,\Omega^1_{\Ss^W})\,|\,
\nabla^{\widehat{J}}\omega=0\}.
\end{equation}
by \EQUref{4.1.3002}. 
Since $\varphi$ is faithfully flat, the sequence
\begin{equation}
0 \to \C \to \Ss^W \to 
\Omega^1_{\Ss^W} \to 
\Omega^2_{\Ss^W} \to 
\cdots\EQUlabel{4.1.5001}
\end{equation}
is exact. Thus we have an exact sequence
\begin{equation}
0 \to \C \to S^W \to 
\Omega^1_{S^W} \to 
\Omega^2_{S^W} \to 
\cdots\EQUlabel{4.1.5000}
\end{equation}
because each homogeneous part of each 
graded module of \EQUref{4.1.5001} is coherent 
and the domain $\Hh$ is Stein. 
Then we could take $t^1,\cdots,t^n \in S^W$ satisfying (1)(i). 
We could take $t^1,\cdots,t^n \in S^W$ so that 
$t^1,\cdots,t^n$ are homogeneous of degree 
$\deg t^1\geq \cdots \geq \deg t^n$. 
Since the Jacobian 
$\frac{\p(t^1,\cdots,t^n)}{\p(s^1,\cdots,s^n)}$ 
is not $0$, degrees of $t^i$ must be $d^i$. 
For a proof of (1)(iii)(iv), see \cite{extendedII}. 

We prove (2). 
For the proof of (2)(i), we first list up the set of degree 
of $s^{\alpha}$.

Put $\{d^1,\cdots,d^n\}=\{p^1,\cdots,p^m\}$ such that 
$1=p^1>p^2>\cdots>p^m=0$. 
We put $Q^i=\{\alpha\,|\,d^{\alpha}=p^i\}$. 

We show $s^{\alpha} \in \Oo(\Hh)[t^1,\cdots,t^{n-1}]$ 
for $\alpha \in Q^i$ by induction on $i$, that is, 
we show it in the order of $i=m$, $i=m-1$, $i=m-2,\cdots$ inductively.

If $i=m$, then $Q^m=\{n\}$ and we have $s^n \in \Oo(\Hh)$, thus 
the assertion is proved for this case.

If $i=m-1$, then for $\alpha \in Q^{m-1}$, we have 
$t^{\alpha}=\sum_{\beta \in Q^{m-1}} f_{\alpha\,\beta}s^{\beta}$ 
with $f_{\alpha\,\beta} \in \Oo(\Hh)$. 
The matrix $(f_{\alpha\,\beta})$ of size $\#Q^{m-1}$ is invertible 
because the Jacobian $(\frac{\p t^{\alpha}}{\p s^{\beta}})$ of 
size $n$ is upper-triangular and invertible. 
Thus $s^{\beta} \in \Oo(\Hh)[t^1,\cdots,t^{n-1}]$ for 
$\beta \in Q^{m-1}$. 

We assume that $s^{\alpha} \in \Oo(\Hh)[t^1,\cdots,t^{n-1}]$ for 
$\alpha \in Q^{i+1}$ ($1 \leq i \leq m-2$). 

Then by the parallel discussion as above, we could show that 
$s^{\alpha}$ ($\alpha \in Q^{i}$) is a linear combination of
$t^{\alpha}$ ($\alpha \in Q^{i}$) modulo $\Oo(\Hh)$-coefficient 
polynomials $s^{\gamma}$ with $\deg s^{\gamma}<p_i$. 
By the assumption of induction, we have 
$s^{\alpha} \in \Oo(\Hh)[t^1,\cdots,t^{n-1}]$ for 
$\alpha \in Q^{i}$. 
Thus we have (2)(i). 

(2)(ii), (2)(iii), (2)(iv) are direct consequences of (2)(i).
\end{proof}
We call these elements $t^1,\cdots,t^n \in S^W$ with the properties 
of Lemma \ref{4.150}(1) {\it the flat coordinates}. \\

\subsection{Flat pencil}\EQUlabel{4.6.1}
The purpose of Section \ref{4.6.1} is 
to recall the notion of a flat pencil. 
We obtain special properties of the Christoffel symbols 
with respect to flat coordinates by the technique of a flat pencil 
by the parallel discussion of \cite{Dubrovin}. 
They are summarized in Proposition \ref{5.3}. 
They will be used to construct a product in Section 4.3. and 
its potential in Section 4.4. 

First we introduce the rational extensions both 
of a symmetric $S^W$-bilinear form and its Levi-Civita connection. 
Let $K(S^W)$ be the quotient field of the integral domain 
$S^W$. We define $\Omega^1_{K(S^W)}$ and 
$\Der_{K(S^W)}$ by 
\begin{equation}
\Omega^1_{K(S^W)}:=K(S^W) \otimes_{S^W} \Omega^1_{S^W},
\quad
\Der_{K(S^W)}:=K(S^W) \otimes_{S^W} \Der_{S^W}.
\end{equation}
Let $h^*:\Omega^1_{S^W} \times \Omega^1_{S^W} \to S^W$ 
be a symmetric $S^W$-bilinear form 
with $0 \neq \det h^*(ds^{\alpha},ds^{\beta})\in S^W$. 
It induces the $K(S^W)$-linear extension of $h^*$:
\begin{align}
	h^*:\,&\Omega^1_{K(S^W)} \times \Omega^1_{K(S^W)} \to K(S^W),
\end{align}
which is non-degenerate because 
$\det h^*(ds^{\alpha},ds^{\beta})$ is a unit in $K(S^W)$. 
The Levi-Civita connection and its dual:
\begin{align}
	\nabla^{h^*}:\,&\Der_{K(S^W)} \times \Der_{K(S^W)} \to 
	\Der_{K(S^W)},\\
	\nabla^{h^*}:\,&\Der_{K(S^W)} \times \Omega^1_{K(S^W)} \to 
	\Omega^1_{K(S^W)}
\end{align}
are defined and characterized by the metric condition 
$\nabla^{h^*}h^*=0$ and torsion free condition 
$\nabla^{h^*}_{\de}\de'-\nabla^{h^*}_{\de'}\de=[\de,\de']$ 
for $\de,\de' \in \Der_{K(S^W)}$.
We call the $K(S^W)$-bilinear form $h^*$ {\it flat}
if the curvature of $\nabla^{h^*}$ vanishes, i.e. 
\begin{equation}
\nabla^{h^*}_{\de}
\nabla^{h^*}_{\de'}
-
\nabla^{h^*}_{\de'}
\nabla^{h^*}_{\de}
=
\nabla^{h^*}_{[\de,\de']}
\end{equation}
for any $\de,\de' \in \Der_{K(S^W)}$. 

We shall come back to our situation. 
We remind that 
\begin{equation}
I^*:\Omega^1_{S^W}\times \Omega^1_{S^W}\to S^W
\end{equation}
is defined as a global section of \EQUref{2.5.1004}.

The $K(S^W)$-linear extension of 
$I^*$ is non-degenerate and flat 
because $I^*_{\EW}$ is non-degenerate and flat 
on the open dense subset $\EEE \subset \EW$ 
by Proposition \ref{3.1.1500}. 

Taking a global section on $\Hh$ of 
the dual tensor $\widehat{J}^*$ of $\widehat{J}$
in \EQUref{4.16.4}, we have 
\begin{equation}
\widehat{J}^*:\Omega^1_{S^W}\times \Omega^1_{S^W}\to S^W.
\end{equation}
The $K(S^W)$-linear extension of $\widehat{J}^*$ is non-degenerate
and flat by Proposition \ref{4.16.1}.

We denote the Levi-Civita connections 
for $K(S^W)$-linear extensions $I^*$ and $\widehat{J}^*$
by $\nabla^{I^*}$ and $\nabla^{\widehat{J}^*}$ respectively. 

Hereafter we use the flat coordinates $t^1,\cdots,t^n \in S^W$
introduced in Lemma \ref{4.150}(1). 

We fix some notations. 
We simply denote $\frac{\p}{\p t^{\alpha}}$ by $\p_{\alpha}$. 
Thus $\widehat{e}=\p_1$. 
We put
\begin{equation}
\eta^{\alpha \beta}:=\widehat{J}^*(dt^{\alpha},dt^{\beta})
\in \C. 
\end{equation} 
We have $\det (\eta^{\alpha \beta})\neq 0$ 
because the set $\{dt^1,\cdots,dt^{n}\}$ is an
$S^W$-free basis of $\Omega^1_{S^W}$, 
and $\widehat{J}$ is non-degenerate. 
The complex numbers $\eta_{\alpha \beta}$ are determined by the property
\begin{equation}
\eta_{\alpha \beta}\eta^{\beta \gamma}=\de^{\gamma}_{\alpha},
\end{equation}
where we take summation for the same letter. \\

We put 
\begin{equation}\EQUlabel{4.01}
g^{\alpha \beta}:=I^*(dt^{\alpha},dt^{\beta}) \in S^W.
\end{equation}
We put 
	\begin{equation}
		\Gamma^{\alpha \beta}_{\gamma}:=
	I^*(dt^{\alpha},\nabla^{I^*}_{\gamma} dt^{\beta}) \in K(S^W),\ 
	\mbox{ where }\nabla^{I^*}_{\gamma}:=\nabla^{I^*}_{\p_{\gamma}}.
	\end{equation}

\begin{Prop}\EQUlabel{5.2}
Let $t^1,\cdots,t^n$ be the flat coordinates defined as above. \\
(1) $\Gamma^{\alpha \beta}_{\gamma}$ is an element of $S^W$. \\
(2) $g^{\alpha \beta}$ and $
	\Gamma^{\alpha \beta}_{\gamma}$ satisfy 
	\begin{equation}\EQUlabel{5.500}
		\p_1^2(g^{\alpha \beta})=0,\quad
		\p_1^2(\Gamma^{\alpha \beta}_{\gamma})=0.
	\end{equation}\\
(3) $\det (\p_1 g^{\alpha \beta})$ is a unit in 
	$S^W$.
\end{Prop}
\begin{proof}
(1) is a direct consequence of the results of 
\cite{extendedII}. 
We only give the outline. 
By \cite[p.43, (6.7)]{extendedII}, 
$\nabla^{I^*}_{\gamma} dt^{\beta}$ becomes a logarithmic form 
in the sense of \cite{extendedII}. 
Meanwhile $I^*(\omega,\omega')$ is an element of $S^W$ 
for $\omega \in \Omega^1_{S^W}$ and a logarithmic form 
$\omega'$ by \cite[p.38, (5.5.1)]{extendedII}. 
Thus we obtain the assertion of (1).

For the proof of (2), we first check the degrees of 
$(\p_1)^2 g^{\alpha \beta}$
and 
$(\p_1)^2 \Gamma^{\alpha \beta}_{\gamma}$. 
We have
	\begin{equation}
	\deg (\p_1)^2 g^{\alpha \beta}
	=d^{\alpha}+d^{\beta}-2 \leq 0, \quad
	\deg (\p_1)^2 \Gamma^{\alpha \beta}_{\gamma}
	=d^{\alpha}+d^{\beta}-d^{\gamma}-2 \leq 0.
	\end{equation}
Their degrees are $0$ only when 
$\alpha=\beta=1,\gamma=n$. 
In this case, 
$(\p_1)^2 g^{11}=(\widehat{e})^2 g^{11}=\widehat{e} \eta^{11}=0$. 
We show $(\p_1)^2 \Gamma^{11}_n=0$. 
Since 
$$
	\Gamma^{11}_n=I^*(dt^1,\nabla^{I^*}_n dt^1)
	=\frac{1}{2}\p_n I^*(dt^1,dt^1)
	=\frac{1}{2}\p_n g^{11},
$$
it follows that 
$
	(\p_1)^2 \Gamma^{11}_n
	=\frac{1}{2}\p_n (\p_1)^2 g^{11}=0
$. 
In the case where degrees are negative, then 
$
	(\p_1)^2 g^{\alpha \beta}
	=(\p_1)^2 \Gamma^{\alpha \beta}_{\gamma}=0
$. 

For a proof of (3), we remind that $\widehat{J}^*$ is non-degenerate. 
Thus $\det J^*(dt^{\alpha},dt^{\beta})=\det(\p_1g^{\alpha \beta})$ 
is a unit in $S^W$.
\end{proof}

We show that 
$I^*$ and $\widehat{J}^*$ give a flat pencil 
in the sense of \cite[p.194, Def. 3.1]{Dubrovin}. 
\begin{Prop}\EQUlabel{5.1}
$K(S^W)$-linear extensions 
\begin{align}
	I^*:\Omega^1_{K(S^W)} \times \Omega^1_{K(S^W)} &\to K(S^W),\\
	\widehat{J}^*:\Omega^1_{K(S^W)} \times \Omega^1_{K(S^W)} &\to K(S^W)
\end{align}
form a flat pencil \cite[p194 (3.35)]{Dubrovin}. Namely, if 
we put $I^*_{\lambda}:=I^*+\lambda \widehat{J}^*$ 
for any $\lambda \in \C$, we have the following.\\ 
(1) $I^*_{\lambda}$
is non-degenerate and flat. \\
(2) Let $\nabla^{\lambda}$ be the Levi-Civita connection 
for $I^*_{\lambda}$. Then the equality
	\begin{equation}
		I^*_{\lambda}(\omega_1,\nabla^{\lambda}_{\de}\omega_2)=
		I^*(\omega_1,\nabla^{I^*}_{\de}\omega_2)+
		\lambda	\widehat{J}^*(\omega_1,\nabla^{\widehat{J}^*}_{\de}\omega_2)
	\end{equation}
	holds for 
	$\omega_1,\omega_2 \in \Omega^1_{S^W}$, 
	$\de \in \Der_{S^W}$. 
\end{Prop}
\begin{proof}
The proof is completely parallel to Lemma D.1 in \cite[p.227]{Dubrovin}. 

We assert that for any $(r,s) \in \C^2 \setminus\{(0,0)\}$, 
the tensor $r g^{\alpha \beta}+s \p_1 g^{\alpha \beta}$ 
is non-degenerate, flat and its Christoffel symbol 
$\Gamma^{\alpha \beta}_{(r,s)\,\gamma}$ equals 
$r \Gamma^{\alpha \beta}_{\gamma}+s\p_1 \Gamma^{\alpha \beta}_{\gamma}$. 

We show that the proposition follows from this assertion. 
We obtain (1) by $(r,s)=(1,\lambda)$ 
because $\p_1 g^{\alpha \beta}=\eta^{\alpha \beta}$. 
If $(r,s)=(0,1)$, then we see that 
$\p_1 \Gamma^{\alpha \beta}_{\gamma}$ is a Christoffel symbol of 
$\p_1 g^{\alpha \beta}=\eta^{\alpha \beta}$. 
Thus we obtain (2) by $(r,s)=(1,\lambda)$. 

We show the assertion. 
Using the flat coordinates, 
we regard $g^{\alpha \beta}$ and $\Gamma^{\alpha \beta}_{\gamma}$
as functions on flat coordinates, i.e. 
\begin{equation}
g^{\alpha \beta}(t^1,\cdots,t^n),\quad
\Gamma^{\alpha \beta}_{\gamma}(t^1,\cdots,t^n).
\end{equation}
\indent
We assume that $r \neq 0$. 
Then 
$rg^{\alpha \beta}(t^1+\frac{s}{r},t^2,\cdots,t^n)$ 
is non-degenerate, flat and its Christoffel symbol is 
$r\Gamma^{\alpha \beta}_{\gamma}(t^1+\frac{s}{r},t^2,\cdots,t^n)$. 
Since $g^{\alpha \beta}(t^1,\cdots,t^n),
\Gamma^{\alpha \beta}_{\gamma}(t^1,\cdots,t^n)$ 
are polynomial functions of degree 1 
with respect to $t^1$ by \EQUref{5.500}, we have
\begin{align}
rg^{\alpha \beta}(t^1+\frac{s}{r},t^2,\cdots,t^n)
&=rg^{\alpha \beta}(t^1,\cdots,t^n)
+s\p_1g^{\alpha \beta}(t^1,\cdots,t^n),\EQUlabel{4.22}\\
r\Gamma^{\alpha \beta}_{\gamma}(t^1+\frac{s}{r},t^2,\cdots,t^n)
&=r\Gamma^{\alpha \beta}_{\gamma}(t^1,\cdots,t^n)
+s\p_1\Gamma^{\alpha \beta}_{\gamma}(t^1,\cdots,t^n).\EQUlabel{4.23}
\end{align}
Thus we proved the assertion for the case $r \neq 0$. 

For the case of $r=0$, $s\p_1g^{\alpha \beta}=s\eta^{\alpha \beta}$
is non-degenerate and flat because $s \neq 0$. 
On Christoffel symbol, we see that 
$\Gamma^{\alpha \beta}_{(r,s)\,\gamma}
-[r \Gamma^{\alpha \beta}_{\gamma}+s\p_1 \Gamma^{\alpha \beta}_{\gamma}]$
is a rational function 
with respect to $(r,s)$. Since it is $0$ on the domain $r \neq 0$, 
we see that it is $0$ for any $(r,s) \in \C^2\setminus \{(0,0)\}$. 
Thus we proved the assertion.
\end{proof}

The following is a direct consequence of Proposition \ref{5.1} 
(cf. \cite[p.226, (D.1a), (D.2)]{Dubrovin}) and \EQUref{4.1.5000}.
\begin{Prop}\EQUlabel{5.3}
(1) There exists a homogeneous element 
	$f^{\beta} \in S^W$ 
	satisfying the following relations
	\begin{equation}
	\Gamma^{\alpha \beta}_{\gamma}
	=\eta^{\alpha \e}\p_{\e}\p_{\gamma}
	f^{\beta}\quad
	(\alpha,\gamma=1,\cdots,n).
	\EQUlabel{5.24}
	\end{equation}
(2) We have
	\begin{equation}
	\Gamma^{\alpha \beta}_{\gamma}
	\Gamma^{\gamma \de}_{\mu}
	=
	\Gamma^{\alpha \de}_{\gamma}
	\Gamma^{\gamma \beta}_{\mu}\quad
	(\alpha,\beta,\de,\mu=1,\cdots,n).
	\end{equation}
\end{Prop}

We use the following results in Section 4.4 and 4.5. 
\begin{Lemma}\EQUlabel{9.5}
We have
\begin{align}
g^{n \alpha}&=
\eta^{1 n}d^{\alpha}t^{\alpha},\EQUlabel{9.011}\\
\Gamma^{\alpha n}_{\beta}&=0\EQUlabel{9.013},\\
\Gamma^{n \alpha}_{\beta}&=
\eta^{1 n}d^{\alpha}\de^{\alpha}_{\beta}.
\EQUlabel{9.012}
\end{align}
\end{Lemma}

\begin{proof}
For (\ref{9.011}), we should prove 
\begin{equation}\EQUlabel{9.500.500}
I^*(dt^n)=\eta^{1 n}E. 
\end{equation}
We define the $S^W$-isomorphism
\begin{equation}
\Der_{S^W}\stackrel{\sim}{\to}\Omega^1_{S^W},\ 
\de \mapsto \widehat{J}(\de,\cdot)
\end{equation}
induced by $\widehat{J}:\Der_{S^W}\times \Der_{S^W} \to S^W$ 
and denote it also by $\widehat{J}$. 
By \cite[p.51, (9.8), Assertion(iii)]{extendedII} and 
\cite[p.52, (9.9), Cor.]{extendedII}, we have 
\begin{equation}
I^*(\widehat{J}(\widehat{e}))=E.
\end{equation}
We remark that the 
Euler field $E$ in 
\cite[p.38, (5.4.3)]{extendedII} corresponds to 
our operator $E'$. 
Then by 
$dt^n=\eta^{kn}\widehat{J}(\frac{\p}{\p t^k})
=\eta^{1n}\widehat{J}(\frac{\p}{\p t^1})
=\eta^{1n}\widehat{J}(\widehat{e})$, we have the result. 

For (\ref{9.013}), we should prove 
\begin{equation}
\nabla^{I^*}_{\beta}dt^n=0
\end{equation}
because $\Gamma^{\alpha n}_{\beta}
=I^*(dt^{\alpha},\nabla^{I^*}_{\beta}dt^n)
$ 
by definition. 
Since we have
$$
\nabla^{I^*}_{\beta}dt^n=\nabla^{I^*}_{\beta}
(I^*)^{-1}(\eta^{1n}E)
=\eta^{1n}(I^*)^{-1}(\nabla^{I^*}_{\beta}E)
=0
$$
by \EQUref{9.500.500} and $\nabla^{I^*}_{\beta}E=0$ 
(cf. \cite[p.43, (6.6)]{extendedII}), we have the result.

For (\ref{9.012}), we have
\begin{equation}
\p_{\beta}g^{n \alpha}=
I^*(\nabla^{I^*}_{\beta}dt^n,dt^{\alpha})+
I^*(dt^n,\nabla^{I^*}_{\beta}dt^{\alpha})
=\Gamma^{\alpha n}_{\beta}+\Gamma^{n \alpha}_{\beta}.
\end{equation}
By (\ref{9.011}) and (\ref{9.013}), we have (\ref{9.012}).
\end{proof}

\subsection{A construction of a product}

The purpose of this subsection is to define 
a product.

Let $t^1,\cdots,t^n$ be the flat coordinates defined right 
after Lemma \ref{4.150}. 

In order to explain the definition of \EQUref{4.600.1}, 
we assume in this paragraph that there exists a product structure 
$\circ$ on the tangent bundle of $\EW$ such that 
$(\EW,\circ,\widehat{e},E,\widehat{J})$ 
becomes a Frobenius manifold. 
Under this assumption, we obtain the following properties.  
We put the structure coefficients 
$C^{\alpha \beta}_{\gamma}$ with respect 
to the $\Oo_{\EW}$-free basis 
$\widehat{J}^*(dt^{1}),\cdots,\widehat{J}^*(dt^{n})$ 
by the equations:
$$
\widehat{J}^*(dt^{\alpha})
\circ
\widehat{J}^*(dt^{\beta})
=C^{\alpha \beta}_{\gamma}
\widehat{J}^*(dt^{\gamma}). 
$$
Then by the uniqueness of the Levi-Civita connection 
with respect to the tensor $I^*_{\EW}$, we have 
(\cite[p.194, Lemma 3.4]{Dubrovin}) 
\begin{equation}\EQUlabel{4.600.11}
\Gamma^{\alpha \beta}_{\gamma}
=
d^{\beta}C^{\alpha \beta}_{\gamma}.
\end{equation}
Also by the equation 
$\widehat{J}^*(dt^{n})=\eta^{1n}\p_1=\eta^{1n}
\widehat{e}$, we have
\begin{equation}\EQUlabel{4.600.12}
C^{\alpha n}_{\gamma}=\eta^{1n}\de^{\alpha}_{\gamma}. 
\end{equation}
By the equations \EQUref{4.600.11}, \EQUref{4.600.12} 
and the fact that $d^{\beta} \neq 0$ if $\beta \neq n$, 
we see that the product structure is unique if it exists.

Therefore we define the new product $\widehat{\circ}$ by 
the equations:
\begin{equation}\EQUlabel{4.600.2}
\widehat{J}^*(dt^{\alpha})\ 
\widehat{\circ}\ 
\widehat{J}^*(dt^{\beta})
:=\widehat{C}^{\alpha \beta}_{\gamma}
\widehat{J}^*(dt^{\gamma})
\end{equation}
where 
\begin{equation}\EQUlabel{4.600.1}
\widehat{C}^{\alpha \beta}_{\gamma}:=
\begin{cases}
\frac{1}{d^{\beta}}
\Gamma^{\alpha \beta}_{\gamma}&\mbox{if $\beta \neq n$},\\
\eta^{1n}\de^{\alpha}_{\gamma}
&\mbox{if $\beta =n$}.
\end{cases}
\end{equation}
This definition does not depend on the choice of the 
flat coordinates $t^1,\cdots,t^n$. 

By $\Ss^W$-linear extension, we have 
\begin{equation}
\widehat{\circ}:\Der_{\Ss^W} \times \Der_{\Ss^W} \to \Der_{\Ss^W}.
\EQUlabel{4.3.2001}
\end{equation} 
By taking a pull-back of \EQUref{4.3.2001} by $\varphi$ 
in \EQUref{3.1.3001}, we define a product 
\begin{equation}
\widehat{\circ}:\Theta_{\EW}\times \Theta_{\EW} \to \Theta_{\EW}.
\end{equation}

We shall show that 
$(\EW,\widehat{\circ},\widehat{e},E,\widehat{J})$ 
becomes a Frobenius manifold 
in the following subsections. 

\subsection{Existence of a potential}

The purpose of this subsection is to show the existence 
of a potential for the product defined in Section 4.3. 
We give it in Proposition \ref{9.1} adding to the ambiguity of 
a potential.

We explain the idea of the construction of a potential $F$. 
We construct a potential of the product by a technique of 
a flat pencil which is similar to the finite Coxeter group case
 \cite{Dubrovin}. 
But our product $\widehat{\circ}$ is defined in a 
case by case manner (cf. \EQUref{4.600.1}). 
Thus we need to check the compatibility conditions 
also in a case by case manner.

Let $t^1,\cdots,t^n$ be the flat coordinates defined right 
after Lemma \ref{4.150}. 

\begin{Prop}\EQUlabel{9.1}
(1) There exists $F \in S^W$ of degree 2 such that 
\begin{equation}
\widehat{J}(X\ \widehat{\circ}\ Y,Z)=
XYZF
\EQUlabel{9.301}
\end{equation}
for flat vector fields $X,Y,Z$ on $\EW$ with respect to $\widehat{J}$. 
Such $F$ is unique up to adding $c (t^1)^2$ for some $c \in \C$. \\
(2) For any $F \in S^W$ satisfying \EQUref{9.301}, we have
\begin{equation}
I^*_{\EW}(\omega,\omega')=
E\widehat{J}^*(\omega)\widehat{J}^*(\omega')F
\EQUlabel{9.300}
\end{equation}
for flat 1-forms $\omega,\omega'$ on $\EW$ with respect to $\widehat{J}$. 
Conversely any degree 2 element $F \in S^W$ 
satisfying \EQUref{9.300} satisfies \EQUref{9.301}. 
\end{Prop}

\begin{proof}

The assertions are all linear with respect to flat 1-forms $\omega,\omega'$ 
and a flat vector field $X,Y,Z$. 
Then we should only prove the following assertions (a),(b),(c): \\
(a) There exists $F \in S^W$ of degree 2 
such that 
\begin{align}
	\widehat{C}^{\alpha \beta}_{\gamma}&=
	\eta^{\alpha \e}\eta^{\beta \mu}
	\p_{\e}\p_{\mu}\p_{\gamma}F\quad 
	(\alpha, \beta, \gamma=1,\cdots,n),
	\EQUlabel{9.55}	\\
	g^{\beta \gamma}
	&=E\eta^{\beta \e}\eta^{\gamma \mu}
	\p_{\e}\p_{\mu}F\quad
	(\beta,\gamma=1,\cdots,n).
	\EQUlabel{9.91}
\end{align}
\noindent
(b) An element $F \in S^W$ satisfying \EQUref{9.55}
is unique up to adding $c(t^1)^2$ for some $c \in \C$.\\
(c) An element $F \in S^W$ satisfying \EQUref{9.91} 
is unique up to adding $c(t^1)^2$ for some $c \in \C$.\\

Here we used notations such as $g^{\alpha \beta}$ etc. 
defined after Lemma \ref{4.150}. 

We prove (a) in five steps.

As the first step, by Proposition \ref{5.3}(1), we could take a 
homogeneous element $f^{\gamma} \in S^W$ 
satisfying the following relations
\begin{equation}
\Gamma^{\alpha \gamma}_{\sigma}
=\eta^{\alpha \e}\p_{\e}\p_{\sigma}
f^{\gamma}\quad
(\alpha,\sigma=1,\cdots,n).
\EQUlabel{9.24}
\end{equation}
We put
\begin{equation}\EQUlabel{4.100}
F^{\gamma}=
\begin{cases}
f^{\gamma}/d^{\gamma} &\mbox{if $\gamma \neq n$},\\
\frac{1}{2}\eta^{1n}\eta_{\alpha \beta}t^{\alpha}t^{\beta}
&\mbox{if $\gamma =n$}.
\end{cases}
\end{equation}
We show that $F^{\gamma}$ satisfies 
\begin{equation}
\widehat{C}^{\alpha \gamma}_{\sigma}
=\eta^{\alpha \e}\p_{\e}\p_{\sigma}
F^{\gamma}\quad
(\alpha,\sigma=1,\cdots,n).
\EQUlabel{9.310}
\end{equation}
If $\gamma=n$, it is O.K. by definition of $C^{\alpha n}_{\sigma}$. 
If $\gamma \neq n$, then it is O.K. by \EQUref{9.24}
and \EQUref{4.100}.
$F^{\gamma} \in S^W$ is homogeneous of degree 
$1+d^{\gamma}=2-(1-d^{\gamma})$.

As the second step, we shall check the equation: 
\begin{equation}
g^{\beta \gamma}=
(d^{\beta}+d^{\gamma})\eta^{\beta \e}\p_{\e}F^{\gamma}.
\EQUlabel{9.26}
\end{equation}
\indent
If $\gamma \neq n$ in (\ref{9.26}), 
then we should only prove the equation 
\begin{equation}
d^{\gamma}g^{\beta \gamma}
=(d^{\beta}+d^{\gamma})
\eta^{\beta \e}\p_{\e}(d^{\gamma}F^{\gamma})
\EQUlabel{9.21}
\end{equation}
because $d^{\gamma} \neq 0$. 
We use the torsion freeness of $\nabla^{I^*}$
(cf. \cite[p.193, (3.27)]{Dubrovin}):
\begin{equation}
g^{\alpha \s}\Gamma^{\beta \gamma}_{\s}
=
g^{\beta \s}\Gamma^{\alpha \gamma}_{\s}.
\EQUlabel{9.22}
\end{equation}
We take $\beta=n$. Then L.H.S. of (\ref{9.22}) becomes 
\begin{equation}
g^{\alpha \s}\Gamma^{n \gamma}_{\s}
=
g^{\alpha \s}(\eta^{1n}d^{\gamma}\de^{\gamma}_{\s})
=
\eta^{1n}d^{\gamma}g^{\alpha \gamma}
\end{equation}
by \EQUref{9.012}. 
R.H.S. of (\ref{9.22}) becomes
\begin{align}
g^{n\s}\Gamma^{\alpha \gamma}_{\s}
&=(\eta^{1n}d^{\s}t^{\s})
(\eta^{\alpha \e}\p_{\e}\p_{\s}f^{\gamma}) 
\mbox{ by (\ref{9.011}) and (\ref{9.24})}\\
&=\eta^{1n}(d^{\s}t^{\s}\p_{\s})
(\eta^{\alpha \e}\p_{\e}f^{\gamma})\nonumber\\
&=\eta^{1n}(d^{\alpha}+d^{\gamma})(\eta^{\alpha \e}\p_{\e}f^{\gamma}) 
\mbox{ by deg $\eta^{\alpha \e}\p_{\e}=-(1-d^{\alpha})$}.\nonumber
\end{align}
Then by $\eta^{1n} \neq 0$, we have (\ref{9.21}).

If $\gamma=n$ in (\ref{9.26}), 
then we should prove the equation
\begin{equation}
g^{\beta n}=
(d^{\beta}+d^{n})\eta^{\beta \e}\p_{\e}F^n.
\EQUlabel{9.23}
\end{equation}
L.H.S. is $\eta^{1n}d^{\beta}t^{\beta}$ by (\ref{9.011}).
R.H.S. is $d^{\beta}\eta^{\beta \e}\p_{\e}F^n$ by $d^n=0$. 
Then (\ref{9.23}) is a consequence of the definition of 
$F^n$.

As the third step, we show that there 
exists a homogeneous element $F \in S^W$
of degree $2d^1=2$ such that 
it satisfies the following equation:
\begin{equation}
F^{\beta}=\eta^{\beta \mu}\p_{\mu}F.
\EQUlabel{9.56}
\end{equation}
By \EQUref{4.1.5000}, we should only prove the integrability conditions
\begin{equation}
\eta^{\beta \e}\p_{\e}F^{\gamma}
=\eta^{\gamma \e}\p_{\e}F^{\beta}
\EQUlabel{9.31}
\end{equation}
for $\beta \neq \gamma$. 
Since $\beta \neq \gamma$, we have $d^{\beta}+d^{\gamma} \neq 0$. 
Then by (\ref{9.26}), the assersion 
(\ref{9.31}) reduces to the property of metric 
$g^{\beta \gamma}=g^{\gamma \beta}$.

As the fourth step, we have \EQUref{9.91} 
because $E\eta^{\beta \e}\eta^{\gamma \mu}
	\p_{\e}\p_{\mu}F=(d^{\beta}+d^{\gamma})\eta^{\beta \e}\eta^{\gamma \mu}
	\p_{\e}\p_{\mu}F$ for $EF=2F$.

As the fifth step, we have \EQUref{9.55} because 
of \EQUref{9.310} and \EQUref{9.56}. 

Thus we finished the proof of the part (a).

We prove the part (b).
Let $F_1$ and $F_2$ be degree 2 elements of $S^W$ 
satisfying the condition \EQUref{9.55}. 
Then $F_3:=F_1-F_2$ satisfies 
$0=\eta^{\alpha \e}\eta^{\beta \mu}
\p_{\e}\p_{\mu}\p_{\gamma}F_3$. 
Thus $F_3$ is a polynomial of 
$t^1,\cdots, t^n$ of degree less than or equal to $2$. 
But by the degree condition, 
$F_3$ must be constant times $(t^1)^2$. 
Thus we see the ambiguity of $F$ satisfying 
the condition \EQUref{9.55}.

We prove the part (c).
Let $F_4$ and $F_5$ be degree 2 elements of $S^W$ 
satisfying the condition \EQUref{9.91}. 
Then $F_6:=F_4-F_5$ satisfies 
$0=E\eta^{\beta \e}\eta^{\gamma \mu}\p_{\e}\p_{\mu}F_6=
(d^{\beta}+d^{\gamma})\eta^{\beta \e}\eta^{\gamma \mu}\p_{\e}\p_{\mu}F_6$, 
where the last equality comes from the degree condition. 
Thus we have 
\begin{equation}
\eta^{\beta \e}\eta^{\gamma \mu}\p_{\e}\p_{\mu}F_6=
\begin{cases}
0 &\mbox{if $(\beta,\gamma) \neq (n,n)$},\\
f&\mbox{if $(\beta,\gamma)=(n,n)$}
\end{cases}\EQUlabel{9.200}
\end{equation}
for some element $f \in S^W$ of degree $0$. 
Thus $F_6=\frac{1}{2}f (t^1)^2+g$, 
where $g \in S^W$ is a linear combination 
of $t^1,\cdots, t^n$ plus constant. 
But by the degree condition, $g$ must be $0$. 
Applying the equation \EQUref{9.200} for the case of 
$(\beta,\gamma)=(1,n)$, we see that $f$ must be a constant. 
Thus we see the ambiguity of $F$ satisfying the condition \EQUref{9.91}. 
\end{proof}

\subsection{Property of the product}
The purpose of this subsection is to 
show the properties of the product.

\begin{Prop}\EQUlabel{9.160.1}
For vector fields $X,Y,Z$ on $\EW$, we have\\
(1) $X\ \widehat{\circ}\ Y=Y\ \widehat{\circ}\ X$. \\
(2) $\widehat{e}\ \widehat{\circ}\ X=X$. \\
(3) $(X\ \widehat{\circ}\ Y)\ \widehat{\circ}\ Z
=X\ \widehat{\circ}\ (Y\ \widehat{\circ}\ Z)$. 
\end{Prop}

\begin{proof}
Let $t^1,\cdots,t^n$ be the flat coordinates defined right 
after Lemma \ref{4.150}.

For (1), it is a direct consequence of 
Proposition \ref{9.1}.

For (2), we need to show 
\begin{equation}
\widehat{J}^*(dt^{\alpha})\ 
\widehat{\circ}\ 
\widehat{e}
=
\widehat{J}^*(dt^{\alpha})
\quad
(\alpha=1,\cdots,n).
\EQUlabel{9.62}
\end{equation}
By definition, we have
$$
\widehat{J}^*(dt^{\alpha})\ 
\widehat{\circ}\ 
\widehat{J}^*(dt^{n})
=
\widehat{C}^{\alpha n}_{\beta}
\widehat{J}^*(dt^{\beta})
=\eta^{1n}\de^{\alpha}_{\beta}\widehat{J}^*(dt^{\beta})
=\eta^{1n}\widehat{J}^*(dt^{\alpha}).
$$
Since $\widehat{J}^*(dt^{n})=\eta^{1n}\widehat{e}$ 
and $\eta^{1n}\neq 0$, we obtain (\ref{9.62}).

For (3), we need to show 
\begin{equation}
	\widehat{C}^{\alpha \beta}_{\gamma}
	\widehat{C}^{\gamma \de}_{\mu}
	=
	\widehat{C}^{\alpha \de}_{\gamma}
	\widehat{C}^{\gamma \beta}_{\mu}
	\quad
	(\alpha,\beta,\de,\mu=1,\cdots,n).
	\EQUlabel{9.61}
\end{equation}
We show \EQUref{9.61}. We have 
	\begin{equation}
	\Gamma^{\alpha \beta}_{\gamma}
	\Gamma^{\gamma \de}_{\mu}
	=
	\Gamma^{\alpha \de}_{\gamma}
	\Gamma^{\gamma \beta}_{\mu}
	\quad
	(\alpha,\beta,\de,\mu=1,\cdots,n)
	\end{equation}
by Proposition \ref{5.3}(2). 

We show
	\begin{equation}\EQUlabel{9.500.1}
	\Gamma^{\alpha \beta}_{\gamma}
	=d^{\beta}\widehat{C}^{\alpha \beta}_{\gamma}.
	\end{equation}
If $\beta \neq 0$, it is O.K. by \EQUref{4.600.1}. 
If $\beta=0$, it is O.K. because both hands are $0$ 
by $d^n=0$ and \EQUref{9.013}. 

By \EQUref{9.500.1}, we have
\begin{equation}
	d^{\beta}d^{\de}
	\widehat{C}^{\alpha \beta}_{\gamma}
	\widehat{C}^{\gamma \de}_{\mu}
	=
	d^{\beta}d^{\de}
	\widehat{C}^{\alpha \de}_{\gamma}
	\widehat{C}^{\gamma \beta}_{\mu}
	\quad
	(\alpha,\beta,\de,\mu=1,\cdots,n).
\end{equation}
Therefore we obtain (\ref{9.61}) for the case 
$d^{\beta}d^{\de}\neq 0$. 

For the case $d^{\beta}d^{\de}=0$, 
the index $\beta$ or $\de$ must be $n$. 
Then we have 
$\widehat{C}^{\alpha n}_{\gamma}=\eta^{1n}\de^{\alpha}_{\gamma}$ 
by definition. 
Then the assertion (\ref{9.61}) is apparent.
\end{proof}

\subsection{Construction of Frobenius manifold structure}
The purpose of this subsection is to construct a 
Frobenius manifold structure.

\begin{Prop}\EQUlabel{4.16.3}
The tuple 
$(\EW,\widehat{\circ},\widehat{e},E,\widehat{J})$ 
is a Frobenius manifold 
satisfying the conditions of Theorem \ref{2.9.1}(1). 
\end{Prop}
\begin{proof}
We shall check the properties of Frobenius manifold. 

We check 
$\widehat{J}(X\ \widehat{\circ}\ Y,Z)
=\widehat{J}(X,Y\ \widehat{\circ}\ Z)$. 
We may assume that $X,Y,Z$ are flat. 
Then 
$
\widehat{J}(X\ \widehat{\circ}\ Y,Z)
=XYZ F
$. 
Also we have 
$\widehat{J}(X,Y\ \widehat{\circ}\ Z)
=\widehat{J}(Y\ \widehat{\circ}\ Z,X)
YZX F$. 
Since $XYZ F=YZX F$, 
we have 
$\widehat{J}(X\ \widehat{\circ}\ Y,Z)
=\widehat{J}(X,Y\ \widehat{\circ}\ Z)$. 

We check that the $(3,1)$-tensor $\widehat{\nabla}\ \widehat{\circ}\ $ 
is symmetric. 
We should only prove 
$\widehat{\nabla}_X(Y\ \widehat{\circ}\ Z)
=\widehat{\nabla}_Y(X\ \widehat{\circ}\ Z)$ 
for flat vector fields $X,Y,Z$. 
We prove 
$\widehat{J}(\widehat{\nabla}_X(Y\ \widehat{\circ}\ Z),W)
=\widehat{J}(\widehat{\nabla}_Y(X\ \widehat{\circ}\ Z),W)$ for a 
flat vector field $W$. 
Since 
$\widehat{J}(\widehat{\nabla}_X(Y\ \widehat{\circ}\ Z),W)
=X\widehat{J}(Y\ \widehat{\circ}\ Z,W)=XYZW F$ 
and 
$\widehat{J}(\widehat{\nabla}_Y(X\ \widehat{\circ}\ Z),W)
=Y\widehat{J}(X\ \widehat{\circ}\ Z,W)=YXZW F$, 
we have the result. 

The flatness of $\widehat{J}$ and the property 
$\widehat{\nabla} \widehat{e}=0$ 
are asserted in Proposition \ref{4.16.1}.

Homogeneity conditions $Lie_E(\ \widehat{\circ}\ )
=1 \cdot\ \widehat{\circ}\ $ and 
$Lie_E(\widehat{J})=\widehat{J}$ (i.e. $D=1$) are consequences of 
$Lie_E F=2$, $Lie_E \widehat{e}=[E,\widehat{e}]=-\widehat{e}$ and 
$Lie_E I^*_{\EW}=0$. 

We prove $I^*_{\EW}(\omega,\omega')
=\widehat{J}(E,\widehat{J}^*(\omega)\ 
\widehat{\circ}\ \widehat{J}^*(\omega'))$. 
We may assume that $\omega,\omega'$ are flat. 
By Proposition \ref{9.1}(2), we have
\begin{equation*}
\widehat{J}(E,\widehat{J}^*(\omega)\ \widehat{\circ}\ 
\widehat{J}^*(\omega'))
=
E\widehat{J}^*(\omega)\widehat{J}^*(\omega')F
=
I^*_{\EW}(\omega,\omega').
\end{equation*}
\end{proof}

\section{Uniqueness of the Frobenius manifold structure}
In this section, we give the proof of Theorem \ref{2.9.1}(2)(3), 
that is,  
the uniqueness of Frobenius manifold structure on $\EW$. 

Theorem\,\ref{2.9.1} (2) is trivial. 
Theorem\,\ref{2.9.1} (3) reduces to Proposition\,\ref{4.2.3}. 
We prepare the following proposition.

\begin{Prop}\EQUlabel{4.2.2}
Let $(M,\circ,e,E,J)$ be a Frobenius manifold 
with intersection form $I^*$. Put 
\begin{align}
\Ff&:=\{\omega \in \Omega^1_M\,|\,Lie_e \omega=0\},\\
T&:=\{f \in \Oo_M\,|\,e(f)=0\},\\
\Omega^{1\nabla}_M&:=\{\omega \in \Omega^1_M\,|\,\nabla \omega=0\},
\end{align}
where $\nabla$ is the Levi-Civita connection for $J$. 
Then \\
(1) $\Ff \supset \Omega^{1\nabla}_M$ and it induces $\Ff\simeq T \otimes_{\C}
\Omega^{1\nabla}_M$. \\
(2) $e$ is non-singular and $[E,e]=-e$. \\
(3) $e^2I^*(\omega,\omega')=0$ for $\omega,\omega' \in \Ff$. \\
(4) $eI^*(\omega,\omega')=J^*(\omega,\omega')
\mbox{ for }
\omega,\omega' \in \Omega_e$. 
\end{Prop}
\begin{proof}
We show (1). 
First we show that $\Ff \supset \Omega^{1\nabla}_M$. 
We take a flat 1-form $\eta \in \Omega^{1\nabla}_{M}$. 
For a flat vector field $Y$, we have 
$(Lie_e \eta)(Y)=e(\eta(Y))-\eta([e,Y])=0$
because $[e,Y]=\nabla_e Y-\nabla_Y e=0$. 
This gives $Lie_{e}\eta=0$. Thus $\eta \in \Ff$. 
We see easily that the isomorphism 
$\Omega^1_M \simeq \Oo_M \otimes_{\C}\Omega^{1\nabla}_M$ 
induces $\Ff \simeq T\otimes_{\C}\Omega^{1\nabla}_M$. 

We show (2). 
Since $e$ is flat, $e$ is non-singular or $0$. 
If $e=0$, then any vector field $X$ must be $0$ 
because $X=X \circ e=X \circ 0=0$, 
which is a contradiction. Thus $e$ is non-singular. 
Also we have $[E,e]=-e$ because the Lie derivative of 
$e \circ e=e$ by $E$ gives $Lie_E(e)=-e$ since $Lie_E(\circ)=1\cdot \circ$. 

We show (3) and (4). 
We first remark that the local existence of $f \in \Oo_{M}$ such that 
\begin{equation}
J(X,Y \circ Z)=XYZf \EQUlabel{5.21}
\end{equation}
for flat fields $X,Y,Z$ is well-known (cf. \cite[p.147]{Hertling}). 

Then for $\omega,\omega'$ : flat 1-forms, we have 
\begin{align*}
eI^*(\omega,\omega')
&=eJ(E,J^*(\omega)\circ J^*(\omega'))\\
&=eE J^*(\omega)J^*(\omega')f\\
&=(e+Ee)J^*(\omega)J^*(\omega')f\\
&=J^*(\omega,\omega')+EJ^*(\omega,\omega')\\
&=J^*(\omega,\omega')
\end{align*}
because $J^*(\omega,\omega')$ is a constant for flat 1-forms 
$\omega,\omega'$. Then we have 
$
e^2I^*(\omega,\omega')
=eJ^*(\omega,\omega')
=0
$. 

By the result of (1), it is sufficient to show (3) and (4) 
only for $\omega,\omega'$ flat 1-forms, because  
(3) and (4) are linear over the ring $T$. 
Thus we have the result. 
\end{proof}

\begin{Prop}\EQUlabel{4.2.3}
Let $(\EW,\circ,e,E,J)$ be any Frobenius manifold 
which satisfies the conditions of Theorem \ref{2.9.1}(1). 
Let $(\EW,\widehat{\circ},\widehat{e},E,\widehat{J})$
 be a Frobenius manifold constructed in Proposition\,\ref{4.16.3}. 
Then there exists $c \in \C^*$ such that 
\begin{equation}
(\EW,c^{-1}\circ,ce,E,c^{-1}J)=
(\EW,\widehat{\circ},\widehat{e},E,\widehat{J}).
\end{equation}
\end{Prop}
\begin{proof}
By Proposition\,\ref{4.2.2}\,(2), we have $e \in \Der^{lowest}_{S^W}$. 
By Proposition\,\ref{4.2.2}\,(3), $e \in V$, 
where $V$ is defined in \EQUref{4.500}. 
By \EQUref{4.2.1}, we have $e=c^{-1}\ \widehat{e}$ for some $c \in \C^*$. 
We have $\Omega_{e}=\Omega_{\widehat{e}}$, 
where $\Omega_e$ is defined in \EQUref{4.501} 
for $e \in \Der^{lowest}_{S^W}$.

By Proposition\,\ref{4.2.2}(4), 
${J}^*(\omega,\omega')=c^{-1}\widehat{J}^*(\omega,\omega')$ 
for $\omega,\omega' \in \Omega_e=\Omega_{\widehat{e}}$. 
Since $\Omega_e=\Ff$ contains an $\Oo_{\EW}$-free basis of 
$\Omega^1_{\EW}$ by Proposition\,\ref{4.2.2}(1), 
we have ${J}^*=c^{-1}\widehat{J}^*$. 
Thus we have $J=c\widehat{J}$. 

By Theorem \ref{2.9.1}(2), 
\begin{equation}
(\EW,\circ',e',E,J'):=
(\EW,c^{-1}{\circ},c{e},E,c^{-1}{J})
\end{equation}
is also a Frobenius manifold satisfying 
the conditions of Theorem \ref{2.9.1}(1). 
We need to prove 
$(\EW,\circ',e',E,J')=
(\EW,\widehat{\circ},\widehat{e},E,\widehat{J})$. 
We already have $e'=\widehat{e}$, $J'=\widehat{J}$. 

Since these Frobenius manifold structures have the common 
intersection form $I^*_{\EW}$, 
the product structure of the Frobenius manifold 
is uniquely determined by the data of 
the unit vector $e$, the Euler field $E$ and the flat metric $J$ 
as we discussed in Section 4.3. 
Therefore we have the result. 
\end{proof}

\end{document}